# On Hypergeometric 3F2(1)

by

Michael Milgram[*]

Box 1484, Deep River, Ont., Canada.  K0J 1P0

**Abstract:** By systematically applying ten inequivalent two-part relations between hypergeometric sums $_3F_2(\ldots|1)$ to the published database of all such sums, 66 new sums are obtained.  Many results extracted from the literature are shown to be special cases of these new sums.  In particular, the general problem of finding elements contiguous to Watson's, Dixon's and Whipple's theorem is reduced to a simple algorithm suitable for machine computation.  Several errors in the literature are corrected or noted.

## 1. Introduction

The evaluation of the hypergeometric sum $_3F_2(\ldots|1)$ is of ongoing interest, since it appears ubiquitously in many physics problems.  Particularly tantalizing is the fact that Gauss' theorem gives a simple result for $_2F_1(\ldots|1)$, and Bühring[1] has classified general limiting cases of $_3F_2(\ldots|1)$ with distinct parametric excess, although Wimp[2] has demonstrated that no representation consisting solely of gamma functions exists for the general case.  With reference to the WZ algorithm[3], one might expect that computer algebra systems would by now have become a repository of existing information about such sums, but attempts to extract evaluations for particular values of $_3F_2(\ldots|1)$ rarely generate a useful outcome (see Section 3), in spite of the availability of Petkovšek's implementation[4] of the WZ approach.  In the case of non-terminating series, where the WZ "certification" procedure does not apply, Koornwinder has shown[5] that the WZ method can be used to obtain 3-part transformation identities, but again it is rare that a computer algebra program will generate a closed form result for a sum of interest. Thus, in practice, the practitioner is usually reduced to scanning published tables or the literature in the hope that a particular problem at hand has a known closed form[6] evaluation.

At the same time, it is recognized that the potential universe of closed-form results is much larger than listed in the standard tables, because of the Thomae relations that couple two $_3F_2(\ldots|1)$ with transformed parameter sets.  Since there are 10 inequivalent forms of the 120 Thomae transformations, (see Appendix A), the potential universe of "knowable" results is 9 times larger than that given in the tables (one of the 10 relations is the identity).  With the advent of computer algebra systems, it has now become possible to explore this universe for potentially new, closed-form results, keeping in mind the dictum[3] that such a database can never be complete.

---

[*] mike@geometrics-unlimited.com

This paper describes such a study. In Section 2 the algorithmic details are given. Section 3 gives general background on specific closed form results that were examined and records a number of errors that appear in the literature (and tables). In Section 4 some observations are made that lead to two new, but conjectured, sums. Appendix A lists the ten inequivalent Thomae relations used and Appendix B lists the newly found results. Appendix C gives some particular results that are useful for obtaining closed forms contiguous to Watson's, Whipple's and Dixon's theorems (see Section 3).

It is important to recognize that the identification of particular results from comparison with a tabulated database is a powerful adjunct to the WZ method, since that method can usually be counted on to generate a recursion formula. Without some means of obtaining the starting values for the recursion, it is of little practical utility. As will be demonstrated later, the conjunction of the two methods gives a new and powerful means of evaluating generalized representations of Dixon's, Watson's and Whipple's theorem using a computer algebra program.

## 2. **Method**

Prudnikov et. al. [7], Section 7.4.4 tabulate many closed form relations for special parametric cases of $_3F_2(…|1)$. Each of the 10 Thomae relations listed in Appendix A was applied to 70[*] of the closed form identities listed there, resulting in 630 possibly new and/or different closed form sums. Each of these results was compared against all the others, by searching for the existence of a valid transformation among all the top and bottom parameters taking account of the symmetry that exists among these parameters. If such a transformation was found, the relation that had the most general parametric set was retained, and the particular (equivalent or special) case was discarded. Additionally, all results having a parametric excess equal to a non-positive integer were discarded, unless that particular case corresponded to a terminating sum, as were all cases (Karlsson-Minton[8,9,10]) where a top parameter exceeded a bottom parameter by a small positive integer. Finally, all terminating sums that reduce to other terminating sums with the same number of terms were removed. This procedure yielded a "base set" of fundamental results obtainable from Ref. 7. The list of input cases was then expanded to include closed form results harvested from the literature for which an acceptable transformation from the existing "base set" parameters could not be found, and the same procedure was applied.

There were 89 survivors of this culling procedure, of which 23 covered the set of evaluations initially input, leaving 66 relations that are thought to be new (but see Section 3). These are listed in Appendix B. All of the entries in Appendix B were tested for novelty using computer algebra simplification commands and for validity using numerical evaluation. In none of these cases was a computer algebra code able to arrive at a closed form sum as given in Appendix B, although in most cases a recursion formula emerged, which was not checked. In a few cases, a code did generate a putative, but incorrect, closed form result – each such candidate was numerically verified to be incorrect for a few values of the parameters.

---

[*] It is recognized that many of these results were special cases of more general parameter sets.

## 3. Input cases

### 3.1 Notation

Throughout, I use the symbols *"n","m"* and *"L"* to represent positive integers giving a strong constraint on the acceptable transformations among parameters. Otherwise all symbols represent arbitrary complex, continuous variables. In many of the formulae quoted, it is recognized that the real part of the parametric excess $(b_1 + b_2 - a_1 - a_2 - a_3)$ of a particular hypergeometric function on the left-hand side must exceed zero in order for an infinite series representation to exist. However, the right-hand side of any such equation is a valid representation of that hypergeometric function by the principle of analytic continuation. Therefore, provided that the parametric excess is not a negative absolute constant, such a limitation is irrelevant and is only specifically quoted as a convenient categorization parameter ("excess").

### 3.2 Comments on Section 7.4.4 of Ref. 7.

Each of the formula copied from Section 7.4.4 was tested numerically to guard against transcription errors. In this way, it was discovered that a number of formulae given in that source are incorrect. In particular, [7.4.4.19] does not satisfy any numerical test for arbitrary values of the parameter "b", except in the case that $b = -n$. Thus this equation was retained (with an obvious change of notation) in the form

$$_3F_2\left(\begin{matrix}a, b, -n \\ b+\tfrac{1}{2}, a-n+\tfrac{1}{2}\end{matrix} \bigg| 1\right) = \frac{\sqrt{\pi}\,\Gamma\!\left(a-n+\tfrac{1}{2}\right)\Gamma\!\left(b+\tfrac{1}{2}\right)\Gamma\!\left(\tfrac{1}{2}+b-a+n\right)}{\Gamma\!\left(a+\tfrac{1}{2}\right)\Gamma\!\left(-n+\tfrac{1}{2}\right)\Gamma\!\left(\tfrac{1}{2}+b-a\right)\Gamma\!\left(\tfrac{1}{2}+b+n\right)} \qquad (1)$$

Equations [7.4.4.38] and [7.4.4.73] do not satisfy any non-trivial numerical tests, and lacking a source reference from which a corrected form might have been obtained, these two results were omitted from consideration.

Several misprints were also discovered. Notably, equation [7.4.4.43] should read

$$_3F_2\left(\begin{matrix}1, 2, a \\ 3, b\end{matrix}\bigg| 1\right) = -2(b-2) + 2(b-1)^2 \psi'(b), \qquad (2)$$

equation [7.4.4.55] should read

$$_3F_2\left(\begin{matrix}a, a-n, a-n \\ a-n+1, a-n+1\end{matrix}\bigg| 1\right) = \frac{\pi(n-1)!(n-a)}{(1-a)_{n-1}\sin a\pi}, \qquad (3)$$

equation [7.4.4.67] becomes

$$_3F_2\left(\begin{array}{c}2,a,a\\a+2,a+2\end{array}/1\right) = -a^2(a+1)^2[(2a-1)\psi'(a)-2] \tag{4}$$

and equation [7.4.4.71] should be

$$_3F_2\left(\begin{array}{c}3,a,a+1\\a+2,a+3\end{array}/1\right) = (1+a)^2(2+a)(-2a^2(a-1)\psi(1,a)+a(2a-1))/4. \tag{5}$$

Based on [Ref. 7, 7.4.4.25] relating a particular $_3F_2(1)$ to a special $_2F_1(-1)$ three further results were developed. These are referred to as *"Prudnikov 7.4.4.25 variation 1-3"* in Appendix B and are respectively given below:

$$_3F_2\left(\begin{array}{c}a,a+\tfrac{1}{2},b\\b+\tfrac{3}{2}-a,b-a+1\end{array}/1\right) = \frac{2^{-b}\pi^{\tfrac{1}{2}}\Gamma(2+2b-2a)\Gamma(2+b-4a)}{\Gamma(2+2b-4a)(2a-1)}$$
$$\times\left(\frac{1}{\Gamma(\tfrac{1}{2}b)\Gamma(\tfrac{3}{2}+\tfrac{1}{2}b-2a)}-\frac{1}{\Gamma(\tfrac{1}{2}b+\tfrac{1}{2})\Gamma(1+\tfrac{1}{2}b-2a)}\right) \tag{6}$$

$$_3F_2\left(\begin{array}{c}a,a+\tfrac{1}{2},1\\\tfrac{3}{2}-a-\tfrac{1}{2}n,\,2-a-\tfrac{1}{2}n\end{array}/1\right) = \frac{\Gamma(3-n-2a)\Gamma(2-n-4a)}{\Gamma(3-n-4a)}$$
$$\times\left(\frac{2^{-4a-n}\Gamma(1-2a)}{\Gamma(2-n-4a)}+\frac{1}{2\,\Gamma(2a)}\sum_{k=0}^{n-1}\frac{(-1)^k\Gamma(2a+k)}{\Gamma(2+k-n-2a)}\right) \tag{7}$$

$$_3F_2\left(\begin{array}{c}a,a+\tfrac{1}{2},1\\\tfrac{3}{2}-a+\tfrac{1}{2}n,\,1-a+\tfrac{1}{2}n\end{array}/1\right) = \Gamma(1-2a)\frac{(1+n-2a)}{(1+n-4a)}$$
$$\times\left(\frac{2^{n-1-4a}\Gamma(1+n-2a)}{\Gamma(1+n-4a)}-\frac{1}{2\,\Gamma(2a-n)}\sum_{k=1}^{n-1}\frac{(-1)^k\Gamma(2a-n+k)}{\Gamma(1+k-2a)}\right) \tag{8}$$

### 3.3 Other results extracted from the literature

Two general results given elsewhere[11] survived the tests described in the previous section. They are quoted below and referred to respectively as *"Ref. 11, Lemma 2.1"* and *"Ref. 11, Lemma 2.2"*:

$$_3F_2\binom{1,n+1,a}{n+2,b}|1) = \frac{(n+1)\Gamma(b)}{\Gamma(a)}\left(\frac{\Gamma(a-n-1)(\psi(b-n-1)-\psi(b-a))}{\Gamma(b-n-1)}\right.$$
$$\left.-\sum_{l=0}^{n-1}\frac{\Gamma(a+l-n)}{\Gamma(b+l-n)(l+1)}\right). \qquad (9)$$

and

$$_3F_2\binom{a,b,c}{n+b,c+1,}|1) = \frac{(b)_n \Gamma(c+1)\Gamma(1-a)}{(b-c)_n \Gamma(c+1-a)} +$$
$$c\Gamma(b+n)\Gamma(c-b+1-n)\sum_{l=0}^{n-1}\frac{\Gamma(n-l-a)(-1)^l}{\Gamma(b+n-a-l)\Gamma(n-l)\Gamma(c-b-n+2+l)}. \qquad (10)$$

Notice that a result of Rao et al.[12] citing Ramanujan's notebook corresponds to the special case $b = c$, $n \to n+1$ of **(10)**[13].

Several years ago, Sharma[14] claimed to have found two new closed forms for a particular $_3F_2(\ldots|1)$, based on special choices of parameters for a new $_4F_3(\ldots|1)$ obtained by evaluating a double series. Unfortunately, the new $_4F_3(\ldots|1)$ given in Sharma's paper does not satisfy any numerical tests, and is clearly incorrect, since the right hand-side is symmetric in (Sharma's) variables $\alpha$ and $\rho$, while the left-hand side is not. Using the transformation algorithm described previously, it is noted that the left-hand side of Sharma's equation (7) is a special case of [Ref. 7, 7.4.4.25] and the left-hand side of Sharma's equation (8) corresponds to a special case of [Ref. 7, 7.4.4.47], which in turn is a special case of Watson's theorem $W_{0,0}$ (see Section 3.4). However, the right-hand sides of Sharma's equations do not correspond to the respective right-hand sides of the quoted results of Ref. 7. For these reasons, Sharma's results were omitted from consideration.

Recently, Exton[15] gave a "new" result for a special case of $_3F_2(\ldots|1)$ based on a purportedly "new" two-term transformation[16] between different $_3F_2(\ldots|1)$. It is easily seen that Exton's "new" transformation, is actually a symmetric permutation of **T3** of Appendix A, and that the "new" result is incorrect[17] by calculating the different limits *a, b* or *c=0* on both sides of Exton's equation (13). For this reason, that result was also omitted. Exton has also given some new identities[18] and sums for $_3F_2(\ldots|1)$ by evaluating one part of a double series. Each of these results [Ref. 18, 1.13 and 1.14] corresponds to a special case of results given in Appendix B. Additionally, the case [3.3, q=2] of another of Exton's new results[19] is a special case of **(B.1)**.

All of the relevant results given in Table 2 of Krupnikov and Kölbig[20] are special cases of results given in Appendix B; however two new results found in a paper of Gessel and Stanton[21] – equations (1.6) and (1.9) respectively - are listed below:

$$_3F_2\left(\begin{matrix}2a, 1-a, -n\\ 2a+2, -a-\frac{1}{2}-\frac{3}{2}n\end{matrix}\bigg|1\right) = \frac{((n+3)/2)_n(n+1)(2a+1)}{(1+\frac{1}{2}(n+2a+1))_n(2a+n+1)} \tag{11}$$

$$_3F_2\left(\begin{matrix}-sa+s+1, a-1, -n\\ a+1, -s(n+a)-n\end{matrix}\bigg|1\right) = \frac{a(1+s+sn)_n(n+1)}{(1+s(a+n))_n(a+n)} \tag{12}$$

The pedigree of **(12)** is interesting – Gessel and Stanton[*] only give a proof for *n>0*, citing an unpublished letter from Gosper to Askey for the case *n<0*. Since **(12)** satisfies numerical tests for *n<-1*, a second version of that result with $n \to -n-1$ was included in the database. See also Section 4.

In addition, Ref. 21 also gives, in equation (5.16), a special case of a $_4F_3$, citing a problem posed by Fields and Luke. Since this equation is of Minton type (an upper parameter exceeds a lower parameter by unity), it can be reduced as follows:

$$_4F_3\left(\begin{matrix}1, 1+a/(s-1), bs-a, -n\\ a/(s-1), b+1, 1-a-ns\end{matrix}\bigg|1\right) = {}_3F_2\left(\begin{matrix}1, bs-a, -n\\ b+1, 1-a-ns\end{matrix}\bigg|1\right)$$
$$- \frac{(s-1)(sb-a)n}{a(b+1)(1-a-ns)} {}_3F_2\left(\begin{matrix}2, 1+bs-a, 1-n\\ b+2, 2-a-ns\end{matrix}\bigg|1\right) \tag{13}$$

The first of the two terms on the right can be summed (Ref. 7, 7.4.4.13) leading to:

$$_3F_2\left(\begin{matrix}2, 1+sb-a, 1-n\\ b+2, 2-a-ns\end{matrix}\bigg|1\right) = \frac{(b+1)(a-1+ns)(a+ns)ab\Gamma(n)\Gamma(1-ns+n+b-sb)}{(sb-a)(s-1)\Gamma(b+n+1)\Gamma(1-a-ns+n)}$$
$$\times \sum_{L=0}^{n} \frac{\Gamma(-a-ns+L)\Gamma(b+L)}{\Gamma(L+1)\Gamma(1-ns+b-sb+L)} \tag{14}$$
$$+ \frac{b(a+ns)(-1+a+ns)(b+1)}{n(s-1)(sb-a)(b+n)}$$

Although this result does not appear to be very profound – it simply transforms a finite sum – when **T1-T9** are applied to this equation a number of new results emerge that reduce infinite series of Minton type into finite sums. These are listed in Appendix B; the genesis of each result is labeled *"Gessel and Stanton, (Eq. 5.16)"*. Notice that **B.38** based on **(14)** generalizes **B.53**, which is based on **(12)**. As well, **T6** applied to **(14)** reduces to **(11)** as a special case.

---

[*] In Ref. 21, Eq. 5.16, which sums a terminating $_4F_3$ tested numerically true for continuous "*n*". As well, "strange" Eqs. 3.11 through 3.16 all reduce to a $_2F_1$ so they were omitted.

### **3.4 <u>Generalized (Contiguous) Watson, Dixon and Whipple's theorem</u>**

Lavoie[22] has also given two "probably new" summation formulae contiguous to Watson's theorem. These were checked with the database; it was discovered that Lavoie's equation (2) is new, and it was therefore added; Lavoie's "new" equation (1) can be obtained from equation (2) of that same paper by the application of Thomae relation **T1** so it was excluded. All other relations quoted in that paper are special cases of results already included in the database[*] as are all the $_3F_2(\ldots|1)$ relations recorded in another of Lavoie's works[23]. For completeness' sake, Lavoie's equation (2) is reproduced below:

$$_3F_2\left(\begin{array}{c}a,b,c\\ \tfrac{1}{2}(a+b+1), 2c+1\end{array}\middle|1\right) = \frac{2^{a+b}\Gamma(\tfrac{1}{2}a+\tfrac{1}{2}b+\tfrac{1}{2})\Gamma(c+\tfrac{1}{2})\Gamma(c-\tfrac{1}{2}a-\tfrac{1}{2}b+\tfrac{1}{2})}{\Gamma(\tfrac{1}{2})\Gamma(a+1)\Gamma(b+1)}$$
$$\times\left(\frac{\Gamma(\tfrac{1}{2}a+1)\Gamma(\tfrac{1}{2}b+1)}{\Gamma(c-\tfrac{1}{2}a+\tfrac{1}{2})\Gamma(c-\tfrac{1}{2}b+\tfrac{1}{2})} - \frac{ab\Gamma(\tfrac{1}{2}a+\tfrac{1}{2})\Gamma(\tfrac{1}{2}b+\tfrac{1}{2})}{4\Gamma(c-\tfrac{1}{2}a+1)\Gamma(c-\tfrac{1}{2}b+1)}\right) \quad \textbf{(15)}$$

Additionally, in several other works Lavoie et. al. give new results for cases contiguous to Dixon's,[24] Watson's[25] and Whipple's[26] theorem. As shown by Bailey[27] in the case of a $_3F_2$, Watson's theorem can be easily obtained by the application of (a symmetric permutation of) **T2** of Appendix A to Dixon's theorem, and Whipple's theorem can be subsequently obtained by the application of **T2** to Watson's theorem. It turns out that these three theorems are similarly closed under the application of any of the other Thomae relations, as are the contiguity identities for Whipple's, Watson's and Dixon's theorems studied by Lavoie et. al. Since the results given by Lavoie et. al. for Whipple's and Dixon's theorems are limited to a small subset of near-diagonal contiguous cases[†], whereas Lewanowicz[28] has independently given a general result valid for all off-diagonal instances of the generalized Watson's theorem, in this sense Watson's theorem may be more fundamental than the other two. Therefore, in principle, only Lewanowicz' general result for contiguous elements of the generalized Watson's theorem needs to be retained in any database, with the other relations obtained by the application of any of the Thomae transformations (see below).

However, the complexity of these results and the embodiment of the WZ method in computer algebra codes leads to another practical method of obtaining any particular element of the generalized result for any of these cases. Consider the elements $X_{m,n}$ contiguous to the (well-poised) Dixon's element ($m=n=0$):

$$X(a,b,c,m,n) = X_{m,n} \equiv {}_3F_2\left(\begin{array}{c}a,b,c\\ 1+m+a-b, 1+m+n+a-c\end{array}\middle|1\right) \quad \textbf{(16)}$$

with

---

[*] Is it possible that Lavoie's results are the source for several entries in Prudnikov's table?
[†] 38 cases for Whipple's theorem, 39 for Dixon's theorem.

$$X(a,b,c,n,m) = X(a,c,b,m+n,-n) \qquad \textbf{(16b)}$$

Computer algebra codes will yield (very lengthy) three part recursion formulae in the variables *m* and *n*. With knowledge of four values of $X_{m,n}$ and using the recursion so-obtained[*], it is then possible to obtain any of the sums contiguous to Dixon's theorem ($X_{0,0}$) for any value of *m* and *n*. The case $X_{0,0}$ is well known and given by [Ref. 7, 7.4.4.21], $X_{0,1}$ is given by **(B.50)**, $X_{0,-1}$ corresponds to [Ref. 7, 7.7.4.20], $X_{1,0}$ is given by [Ref. 7, 7.7.4.22] and $X_{1,-1}$ is given (symmetrically) by **(B.50)**.

A similar method can be used to obtain generalizations of Watson's theorem, where the recursion is far simpler.  It is found that any of the elements $W_{m,n}$ defined by the generalized Watson theorem

$$W_{m,n} \equiv {}_3F_2 \left( \begin{matrix} a,b,c \\ \tfrac{1}{2}(a+b+1+m), 2c+n \end{matrix} \middle| 1 \right) \qquad \textbf{(17)}$$

obey three-part recursion formulae in "*m*" and "*n*" given by **(C.1)** and **(C.2)** so it remains necessary to locate eight closed forms to start the recursion for all values of "*m*" and "*n*". The element $W_{0,0}$ corresponding to Watson's theorem is well-known and given by [Ref. 7, 7.4.4.18] or [Ref. 30, 3.13.3(7)]; $W_{0,1}, W_{0,-1}, W_{-1,1}, W_{1,0}$ and $W_{1,-1}$ are given by **(15)**, **(B.43), (B.44), (B.47) and (B.51)** respectively.  Other starting results for the initial elements $X_{m,n}$ can be obtained from Lavoie et. al.[24], and corresponding elements $W_{m,n}$ follow from **(21)**. As a test, the elements $X_{-1,0}$ and $X_{2,0}$, obtained by recursion (listed in **(C.3)** and **(C.4)**), were compared to, and agreed with, the same elements taken from Ref. 24. Using **(21)** (below) the corresponding elements $W_{-1,0}$ and $W_{2,0}$ were found; these are given by **(C.5)** and **(C.6)** respectively.  Using this method any of the other elements $X_{m,n}$ or $W_{m,n}$ can now be found by recursion in either $X_{m,n}$ or (preferably) $W_{m,n}$ using **(C.1)** and **(C.2)**.

The third member of the trio consists of generalized elements contiguous to Whipples' theorem. As noted by Lavoie et. al.[29] using a method similar to that employed here, elements $P_{m,n}$ of Whipple's theorem

$$P_{m,n} \equiv {}_3F_2 \left( \begin{matrix} a, \ b, 1-b+m+n \\ c, 1+2a+m-c \end{matrix} \middle| 1 \right) \qquad \textbf{(18)}$$

---

[*] Note that Lewanowicz' coefficients are also given recursively.

are related to $X_{m,n}$. The relationship can be found by applying **T8** to $X_{m,n}$ giving

$$P_{m,n} = \frac{X(2a-b-c+m+1, 1+a-c+m, 1-b+m+n, m, n)\Gamma(a-n)\Gamma(c)}{\Gamma(b+c-1-m-n)\Gamma(a-b+m+1)}. \tag{19}$$

The correspondence with Watson's theorem is obtained by applying **T1** to $W(a,b,c,m,n) \equiv W_{m,n}$ giving

$$P_{m,n} = \frac{W(c-b, 2a-c+m+1-b, a-n, m, n)\Gamma(a-n)\Gamma(c)\Gamma(1+m+2a-c)}{\Gamma(b)\Gamma(2a-n)\Gamma(a-b+m+1)}. \tag{20}$$

Notice that **(19)** fails[30] for the case $m=n=0$ when "$a$" is a non-positive integer, unless "$b$" is an integer*.

As an adjunct to the above, apply **T2** to $X_{m,n}$ to find

$$W_{m,n} = X(2c+n-a, (b-a+m+1)/2, (1+m-a-b)/2+c+n, m, n)$$
$$\times \frac{\Gamma(c+n+\tfrac{1}{2}(1+m-a-b))\Gamma(2c+n)\Gamma(\tfrac{1}{2}(a+b+1+m))}{\Gamma(a)\Gamma(c+n+\tfrac{1}{2}(1+m+b-a))\Gamma(2c+n+\tfrac{1}{2}(1+m-a-b))}, \tag{21}$$

or

$$X_{m,n} = W(1+m+a-2b, a, 1+m+a-b-c, m, n)$$
$$\times \frac{\Gamma(a-2b-2c+2+2m+n)\Gamma(1+a-c+m+n)}{\Gamma(1-c+m+n)\Gamma(2a-2b-2c+2+2m+n)}. \tag{22}$$

See also **(B.48)** and **(B.49)**.

### 3.5 Other computer applications

Although others have also made extensive use of computer algebra to obtain transformations and identities between generalized hypergeometric series, I found no other new sums in the relevant database. In particular, Gessel[31] gives a large number of new identities of which only one (equation 30.1) is relevant to this work. That result is of the Karlsson-Minton type being a special case of [Ref. 7, 7.4.4.10].

Krattenthaler's Mathematica package HYP[32] and Gauthier's Maple package HYPERG,[33] contain a number of rules taken from Slater[34]. Although not obvious, several of these are special cases of Appendix B. The correspondence follows:

---

* see ref. 30 Eqs 3.13(8) and (9). This exceptional case is given by Prudnikov 7.4.4.104.

**Table 1. Showing the correspondence between the database sums of refs. 32, 33, 34 and the database elements described herein.**

| Equation reference | database number and equivalent special case |
|---|---|
| S3201 | **B.3** with $m=1$ |
| S3204 | [Ref. 7, 7.4.4.10] with $e = 1 - b + \lambda$; $c = 1 + \lambda/2$ |
| S3235 | **B.3** with $m = 2$; $c = 1 + a - b$ |
| Slater III.15 | [Ref. 7, 7.4.4.10] with $c = 1 + \lambda/2$; $b = -n$ |

## 4. Comments and Conjectures

As discussed, Appendix B gives new summation results for a number of hypergeometric sums. Although based on the two-part Thomae relations, many of the results that contain a series are also equivalent to the well-known three part relations between hypergeometric $_3F_2(1)$[30]. However, all these cases still reduce an infinite series with a parameter *"n"* to a finite series with *"n"* (or less) terms and so remain useful, particularly for special values of *"n"*. In a number of instances, new results were obtained which were deleted from Appendix B because they trivially reduced to Karlsson-Minton type; exceptions are **(B.42), (B.53)** and **(B.60-66)** because of the relatively complicated nature of the relations among their parameters.

In a recent work, Krattenthaler and Rivoal[35] introduced two new three-parameter two-part relations between $_3F_2(1)$. These were tested against the database, and 31 results were found that satisfied one or the other of these relations. However, when the appropriate equality was applied to any of these cases in the same manner as was done for the Thomae relations, the resulting sum was found to be already included in the database, so no new sums were found from this procedure.

In the course of testing the results, an interesting observation was made that leads to a conjecture about two possibly new sums. It was discovered that **(11)**, due to Gessel and Stanton[21], satisfies numerical tests for $n \to -2n$, leading to the following conjecture:

$$_3F_2\left({2a, 1-a, 2n \atop 2a+2, -a+3n-\tfrac{1}{2}} \bigg| 1\right) = \frac{(2a+1)\Gamma(a-n+\tfrac{1}{2})\Gamma(-3n+\tfrac{5}{2})}{3\Gamma(\tfrac{3}{2}-3n+a)\Gamma(\tfrac{3}{2}-n)} \quad n \geq 0 \tag{23}$$

which would be a new result. If $n = 1$, **(23)** reduces to **(B.56)** (with the substitutions $a \to 2a$, $b \to 1-a$).

It was also observed that [Ref. 21, (5.16)] survived all numerical tests when *"n"* was not integral. Setting $a = s - 1$ in that equation reduces it to a $_3F_2$, which turns out to be equivalent to **(B.54)**. Replace *"n"* by a continuous variable and with a change of notation arrive at the following conjecture that satisfies all numerical tests that were tried:

$$_3F_2\left(\begin{matrix}2, a, b \\ c, (2c-3-a+ab-b)/(c-2)\end{matrix}\bigg|1\right) = -\frac{(c-1)(c+ab-b-a-1)}{(c-1-b)(a+1-c)}. \tag{24}$$

**(24)** holds true if $b = c$ as well as all other limits that reduce the left-hand side to a $_2F_1$.

### 5. Summary

Sixty-six new sums for hypergeometric $_3F_2(\ldots|1)$ were found. Together with twenty-three previously known forms, this collection forms a useful basic database for the possible computerized identification of any desired sum of the form $_3F_2(\ldots|1)$. This can be achieved by simply seeking a transformation between the parameters of the candidate and any of the elements of the database. The existence of such a database is an important adjunct to the WZ recursion method of evaluating such sums, since it provides a starting point for the recursion in several important cases. Two conjectured new forms were found and a number of errors appearing in the literature were noted.

**Acknowledgement:** This work was sponsored by no one. I am grateful to R. Israel for running verification tests using the Maple code.

**Appendices follow on next 25 pages**

## Appendix A

The following lists 10 inequivalent Thomae relations for $_3F_2([a, b, c], [f, e], 1)$ as they are numbered and referenced in the text. There are 110 other symmetric permutations. T10 is the identity.

T1:

$$_3F_2([e + f - a - b - c, f - c, e - c], [e - b + f - c, e + f - a - c], 1)$$
$$\times \Gamma(e + f - a - b - c)\Gamma(f)\Gamma(e)/(\Gamma(c)\Gamma(e - b + f - c)\Gamma(e + f - a - c))$$

T2:

$$_3F_2([e + f - a - b - c, -b + f, e - b], [e - b + f - c, e - b + f - a], 1)$$
$$\times \Gamma(e + f - a - b - c)\Gamma(f)\Gamma(e)/(\Gamma(b)\Gamma(e - b + f - c)\Gamma(e - b + f - a))$$

T3:

$$\frac{_3F_2([f - c, -b + f, a], [e - b + f - c, f], 1)\,\Gamma(e + f - a - b - c)\,\Gamma(e)}{\Gamma(e - a)\,\Gamma(e - b + f - c)}$$

T4:

$$\frac{_3F_2([e - c, e - b, a], [e - b + f - c, e], 1)\,\Gamma(e + f - a - b - c)\,\Gamma(f)}{\Gamma(f - a)\,\Gamma(e - b + f - c)}$$

T5:

$$_3F_2([e + f - a - b - c, f - a, e - a], [e + f - a - c, e - b + f - a], 1)$$
$$\times \Gamma(e + f - a - b - c)\Gamma(f)\Gamma(e)/(\Gamma(a)\Gamma(e + f - a - c)\Gamma(e - b + f - a))$$

T6:

$$\frac{_3F_2([f - c, f - a, b], [e + f - a - c, f], 1)\,\Gamma(e + f - a - b - c)\,\Gamma(e)}{\Gamma(e - b)\,\Gamma(e + f - a - c)}$$

T7:

$$\frac{_3F_2([e - c, e - a, b], [e + f - a - c, e], 1)\,\Gamma(e + f - a - b - c)\,\Gamma(f)}{\Gamma(-b + f)\,\Gamma(e + f - a - c)}$$

T8:

$$\frac{_3F_2([-b + f, f - a, c], [e - b + f - a, f], 1)\,\Gamma(e + f - a - b - c)\,\Gamma(e)}{\Gamma(e - c)\,\Gamma(e - b + f - a)}$$

T9:

$$\frac{{}_3\mathrm{F}_2([e-b,\,e-a,\,c],\,[e-b+f-a,\,e],\,1)\,\Gamma(e+f-a-b-c)\,\Gamma(f)}{\Gamma(f-c)\,\Gamma(e-b+f-a)}$$

T10:

$${}_3\mathrm{F}_2([a,\,b,\,c],\,[f,\,e],\,1)$$

## Appendix B

The following lists 66 new evaluations of ${}_3\mathrm{F}_2([a,b,c],[e,f],1)$ obtained as discussed in the text and labelled B.1 to B.66. Each is classified according to the (parametric) "excess". The derivation can be duplicated by applying the appropriate Thomae relation from Appendix A to the basic relation whose source is given. For example, B.3 can be obtained by applying T7 of Appendix A to Eq. 7.4.14 of Ref. 7. Limitations on the parameters are given, if necessary, with $m,n > 0$ always.

$$\text{excess} = -b+c-a-n$$

**B.1  Prudnikov 7.4.4.14 requires n<m   : T2**

$${}_3\mathrm{F}_2([a,\,m,\,b],\,[c,\,m-n],\,1) =$$

$$\Gamma(a+n-m+1)\left(\sum_{L=0}^{m-1}\frac{(-1)^L\,\Gamma(-b+c-a-n+L)\,\Gamma(1-b+L)}{\Gamma(m-L)\,\Gamma(-b-n+1+L)\,\Gamma(-b+c-m+1+L)\,\Gamma(L+1)}\right)\Gamma(c)$$
$$\Gamma(m-n)\,(-1)^n/(\Gamma(c-a)\,\Gamma(a))$$

$$\text{excess} = -b+c-a+m$$

**B.2  Prudnikov 7.4.4.14 requires n<m   : T5**

$${}_3\mathrm{F}_2([a,\,b,\,-n],\,[c,\,m-n],\,1) =$$

$$\Gamma(a+n-m+1)\left(\sum_{L=0}^{m-1}\frac{(-1)^L\,\Gamma(c-a+L)\,\Gamma(1+b-m+n+L)}{\Gamma(m-L)\,\Gamma(b-m+1+L)\,\Gamma(c+n-m+1+L)\,\Gamma(L+1)}\right)\Gamma(c)$$
$$\Gamma(m-n)\,(-1)^n/(\Gamma(c-a)\,\Gamma(a))$$

$$\text{excess} = m$$

**B.3  Prudnikov 7.4.4.14   : T7**

$$_3F_2([a, -n, b], [c, a-c+b-n+m], 1) =$$

$$\Gamma(c-b+n-m+1)\,\Gamma(-a+c-b-m+1)$$
$$\left(\sum_{L=0}^{m-1} \frac{(-1)^L\,\Gamma(b+L)\,\Gamma(1-a+c-m+n+L)}{\Gamma(m-L)\,\Gamma(-a+c-m+1+L)\,\Gamma(c+n-m+1+L)\,\Gamma(L+1)}\right)\Gamma(m)\,\Gamma(c)/$$
$$(\Gamma(n+1-a+c-b-m)\,\Gamma(b)\,\Gamma(c-b))$$

$$\text{excess} = -b+c-a+n$$

### B.4  Ref. 11, Lemma 2.2  : T1

$$_3F_2([a, 1, b], [n+1, c], 1) =$$

$$\frac{\Gamma(a-n)\,\Gamma(b-n)\,\Gamma(n+1)\,\Gamma(c)\,\Gamma(-b+c-a+n)}{\Gamma(b)\,\Gamma(c-b)\,\Gamma(a)\,\Gamma(c-a)}$$

$$+\frac{\Gamma(-b+1)\left(\sum_{L=0}^{n-1}\frac{\Gamma(-L+a-1)\,(-1)^L}{\Gamma(c-1-L)\,\Gamma(n-L)\,\Gamma(-b+2+L)}\right)\Gamma(n+1)\,\Gamma(c)}{\Gamma(a)}$$

$$\text{excess} = -b+c-a+1$$

### B.5  Ref. 11, Lemma 2.2  : T2

$$_3F_2([a, b, n], [n+1, c], 1) =$$

$$\frac{\Gamma(b-n)\left(\sum_{L=0}^{n-1}\frac{\Gamma(-L+a-1)\,(-1)^L}{\Gamma(-b+c-L)\,\Gamma(n-L)\,\Gamma(b+1+L-n)}\right)\Gamma(n+1)\,\Gamma(c)\,\Gamma(-b+c-a+1)}{\Gamma(a)\,\Gamma(c-a)}$$

$$+\frac{\Gamma(a-n)\,\Gamma(-b+1)\,\Gamma(n+1)\,\Gamma(c)}{\Gamma(n+1-b)\,\Gamma(c-n)\,\Gamma(a)}$$

$$\text{excess} = n$$

### B.6  Ref. 11, Lemma 2.2  : T6

$$_3F_2([1, a, b], [c, -c+b+a+n+1], 1) =$$

$$\frac{(-c+b+a+n)\,\Gamma(-c+a+1)\left(\sum_{L=0}^{n-1}\frac{\Gamma(-L+c-b-1)\,(-1)^L}{\Gamma(c-1-L)\,\Gamma(n-L)\,\Gamma(-c+a+2+L)}\right)\Gamma(n)\,\Gamma(c)}{\Gamma(c-b)}$$

$$+\frac{\Gamma(-c+b+a+n+1)\,\Gamma(-b+c-n)\,\Gamma(c-a-n)\,\Gamma(n)\,\Gamma(c)}{\Gamma(c-a)\,\Gamma(a)\,\Gamma(c-b)\,\Gamma(b)}$$

$$\text{excess} = 1$$

**B.7  Ref. 11, Lemma 2.2  : T9**

$$_3F_2([n, a, b], [c, -c + b + a + n + 1], 1) =$$

$$\frac{\Gamma(-c + b + a + n + 1)\,\Gamma(c - a - n)\left(\sum_{L=0}^{n-1} \frac{\Gamma(-L + c - b - 1)\,(-1)^L}{\Gamma(a - L)\,\Gamma(n - L)\,\Gamma(c + 1 - a - n + L)}\right)\Gamma(c)}{\Gamma(c - b)\,\Gamma(b)}$$

$$+ \frac{\Gamma(-c + b + a + n + 1)\,\Gamma(-b + c - n)\,(-1)^n\,\Gamma(c - a - n)\,\Gamma(c)}{\Gamma(c - a)\,\Gamma(c - n)\,\Gamma(c - b)\,\Gamma(-c + b + 1 + a)}$$

$$\text{excess} = 1$$

**B.8  Prudnikov 7.4.4.31  : T5**

$$_3F_2([a, b, \tfrac{1}{2}a - b - \tfrac{3}{2} - n], [\tfrac{1}{2}a + b + \tfrac{1}{2}, a - b - 1 - n], 1), =$$

$$\frac{1}{2} \frac{\Gamma(\tfrac{1}{2}a + b + \tfrac{1}{2})\,\Gamma(b + n + 2)\,\Gamma(a - b - 1 - n)\,\pi}{\cos(\tfrac{1}{2}\pi(2b - a))\,\Gamma(-\tfrac{1}{2}a + b + \tfrac{1}{2})\,\Gamma(\tfrac{1}{2}a + \tfrac{1}{2})^2\,\Gamma(2 + n + 2b)\,\Gamma(-2b - n + a - 1)}$$

$$+ \frac{\tfrac{1}{2}\left(\sum_{L=0}^{n+1} \frac{\Gamma(-\tfrac{1}{2}a + b + \tfrac{1}{2} + L)\,\Gamma(-b - n - 1 + L)}{\Gamma(\tfrac{1}{2}a - b - \tfrac{1}{2} - n + L)\,\Gamma(b + 1 + L)}\right)\Gamma(\tfrac{1}{2}a + b + \tfrac{1}{2})\,\Gamma(a - b - 1 - n)}{\Gamma(-b - n - 1)\,\Gamma(a)\,\Gamma(-\tfrac{1}{2}a + b + \tfrac{1}{2})}$$

$$\text{excess} = 1$$

**B.9  Prudnikov 7.4.4.32  : T5**

$$_3F_2([a, b, n + \tfrac{1}{2}a - b - \tfrac{3}{2}], [\tfrac{1}{2}a + b + \tfrac{1}{2}, n + a - b - 1], 1), =$$

$$-\tfrac{1}{2}\Gamma(\tfrac{1}{2}a - b + \tfrac{1}{2})\,\Gamma(b - n + 2)\left(\sum_{L=1}^{n-2} \frac{\Gamma(\tfrac{3}{2} - \tfrac{1}{2}a + b - n + L)\,\Gamma(-b + L)}{\Gamma(\tfrac{1}{2} + \tfrac{1}{2}a - b + L)\,\Gamma(2 + b - n + L)}\right)$$

$$\times \Gamma(\tfrac{1}{2}a + b + \tfrac{1}{2})\,\Gamma(n + a - b - 1)\sin(\pi b)\,(-1)^n \cos(\pi(\tfrac{1}{2}a - b))\Big/(\Gamma(a)\,\pi^2)$$

$$+ \frac{\tfrac{1}{2}\,\Gamma(\tfrac{1}{2}a - b + \tfrac{1}{2})\,\Gamma(\tfrac{1}{2}a + b + \tfrac{1}{2})\,\Gamma(b - n + 2)\,\Gamma(n + a - b - 1)}{\Gamma(\tfrac{1}{2}a + \tfrac{1}{2})^2\,\Gamma(n + a - 2b - 1)\,\Gamma(-n + 2b + 2)}$$

excess= 1

### B.10  Prudnikov 7.4.4.25 variation 2   : T1

$$_3F_2([4a+n, a, a+\frac{1}{2}], [\frac{1}{2}+3a+\frac{1}{2}n, 1+3a+\frac{1}{2}n], 1) =$$

$$\frac{2^{(2a)}\,\Gamma(1+3a+\frac{1}{2}n)\,\Gamma(\frac{1}{2}+3a+\frac{1}{2}n)\left(\sum_{k=0}^{n-1}\frac{(-1)^k\,\Gamma(1-2a-n+k)}{\Gamma(1+2a+k)}\right)}{\Gamma(1-2a-n)\,\sqrt{\pi}\,\Gamma(1+4a+n)}$$

$$+\frac{2^{(-1+6a+n)}\,\Gamma(2a+n)\,\Gamma(1+3a+\frac{1}{2}n)\,\Gamma(\frac{1}{2}+3a+\frac{1}{2}n)}{\Gamma(4a+n)\,\sqrt{\pi}\,\Gamma(1+4a+n)}$$

excess= 1

### B.11  Prudnikov 7.4.4.25 variation 3   : T1

$$_3F_2([1-n+4a, a+\frac{1}{2}, a], [1-\frac{1}{2}n+3a, -\frac{1}{2}n+3a+\frac{3}{2}], 1) =$$

$$\frac{(-\frac{1}{2}n+\frac{1}{2}+2a)\,\Gamma(1+2a-n)\,\Gamma(2-n+6a)}{\Gamma(2-n+4a)^2}$$

$$-\frac{\Gamma(1+2a-n)\left(\sum_{k=1}^{n-1}\frac{(-1)^k\,\Gamma(-2a+k)}{\Gamma(1-n+2a+k)}\right)\Gamma(1-\frac{1}{2}n+3a)\,\Gamma(-\frac{1}{2}n+3a+\frac{3}{2})}{\Gamma(-2a)\,\Gamma(2-n+4a)\,\Gamma(a+\frac{1}{2})\,\Gamma(1+a)}$$

excess= $-\frac{1}{2}b-\frac{1}{2}n$

### B.12  Prudnikov 7.4.4.31   : T1

$$_3F_2([a, b, \frac{1}{2}a+\frac{1}{2}], [1+a, \frac{1}{2}a-\frac{1}{2}n-\frac{1}{2}+\frac{1}{2}b], 1) =$$

$$\frac{1}{2}\frac{\Gamma(\frac{1}{2}a-\frac{1}{2}n-\frac{1}{2}+\frac{1}{2}b)\,a\left(\sum_{L=0}^{n+1}\frac{\Gamma(-\frac{1}{2}a-\frac{1}{2}n-\frac{1}{2}+\frac{1}{2}b+L)\,\Gamma(-\frac{1}{2}b-\frac{1}{2}n+L)}{\Gamma(\frac{1}{2}a-\frac{1}{2}n+\frac{1}{2}-\frac{1}{2}b+L)\,\Gamma(\frac{1}{2}b-\frac{1}{2}n+L)}\right)}{\Gamma(-\frac{1}{2}a-\frac{1}{2}n-\frac{1}{2}+\frac{1}{2}b)}$$

$$-\frac{1}{2}\frac{\Gamma(\frac{1}{2}a-\frac{1}{2}n-\frac{1}{2}+\frac{1}{2}b)\,\Gamma(\frac{1}{2}a+\frac{1}{2}n+\frac{3}{2}-\frac{1}{2}b)\,\sqrt{\pi}\,2^a\,\Gamma(1+\frac{1}{2}a)}{\Gamma(b)\,\Gamma(1+a-b)\,\Gamma(\frac{1}{2}a+\frac{1}{2})\sin(\frac{1}{2}\pi(b+n))}$$

$$\text{excess} = b - \frac{1}{2}n - \frac{1}{2}a - 1$$

### B.13  Prudnikov 7.4.4.31  : T3

$$_3F_2([a,\, b,\, 1],\, [2b,\, \frac{1}{2}a - \frac{1}{2}n],\, 1) =$$

$$\Gamma(\frac{1}{2}a - \frac{1}{2}n)\,(b - \frac{1}{2})\,\Gamma(b - \frac{1}{2}n - \frac{1}{2}a - 1)$$

$$\left( \sum_{L=0}^{n+1} \frac{\Gamma(-b - \frac{1}{2}n + \frac{1}{2}a + L)\,\Gamma(-\frac{1}{2}a - \frac{1}{2}n + L)}{\Gamma(b - \frac{1}{2}n - \frac{1}{2}a + L)\,\Gamma(\frac{1}{2}a - \frac{1}{2}n + L)} \right) \Big/ \left( \Gamma(\frac{1}{2}a - b - \frac{1}{2}n)\,\Gamma(-\frac{1}{2}a - \frac{1}{2}n) \right)$$

$$- \frac{\sqrt{\pi}\,2^{(2b-2)}\,\Gamma(\frac{1}{2}a - \frac{1}{2}n)\,\Gamma(b + \frac{1}{2}n - \frac{1}{2}a + 1)\,\Gamma(b + \frac{1}{2})\,\Gamma(b - \frac{1}{2}n - \frac{1}{2}a - 1)}{\Gamma(2b - a)\,\Gamma(a)\,\sin(\frac{1}{2}\pi(a + n))\,\Gamma(-\frac{1}{2}a - \frac{1}{2}n)\,\Gamma(b)}$$

$$\text{excess} = -n - 2b$$

### B.14  Prudnikov 7.4.4.31  : T6

$$_3F_2([a,\, \frac{1}{2}a - \frac{1}{2}n - 1,\, b],\, [a - b - 1 - n,\, \frac{1}{2}a - \frac{1}{2}n],\, 1) =$$

$$\left( \left( -\frac{1}{2} \frac{\left( \sum_{L=0}^{n+1} \frac{\Gamma(b + L)\,\Gamma(-\frac{1}{2}a - \frac{1}{2}n + L)}{\Gamma(-b - n + L)\,\Gamma(\frac{1}{2}a - \frac{1}{2}n + L)} \right) \sin(\frac{1}{2}\pi(a + n))\,\Gamma(\frac{1}{2}a - \frac{1}{2}n)\,\Gamma(-n - 2b)}{\pi\,\Gamma(-2b - n + a - 1)} \right. \right.$$

$$\left. + \frac{\pi^{(3/2)}\,(-1)^n\,\Gamma(-a + 1)}{\sin(\pi b)\,2^{(a-n)}\,\Gamma(\frac{1}{2}a - b - \frac{1}{2}n)^2\,\Gamma(\frac{1}{2} + \frac{1}{2}a - \frac{1}{2}n)\,\Gamma(n + 1 - a)} \right)$$

$$\left. \times \Gamma(\frac{1}{2}a + \frac{1}{2}n + 1)\,\Gamma(a - b - 1 - n) \right) \Big/ \Gamma(b)$$

$$\text{excess} = a$$

### B.15  Prudnikov 7.4.4.31  : T7

$$_3F_2([a,\, b,\, -b - n - 1],\, [-b - n - 2 + 2a,\, b + 1],\, 1) =$$

$$\frac{1}{2}\left(\sum_{L=0}^{n+1}\frac{\Gamma(-b-n-1+L)\,\Gamma(b-a+1+L)}{\Gamma(b+1+L)\,\Gamma(a-b-1-n+L)}\right)$$

$$\times\Gamma(a)\,\Gamma(b+1)\,\Gamma(-b-n-2+2\,a)/(\Gamma(2\,a-1)\,\Gamma(-b-n-1)\,\Gamma(b-a+1))$$

$$+2^{(2\,n-1)}\,\sqrt{\pi}\,\Gamma(-1-n-2\,b)\,\Gamma(3+n-2\,a+2\,b)\,\Gamma(b+n+2)\,\Gamma(a-b)\,\Gamma(a-\tfrac{1}{2})\,\Gamma(-b-n-2+2\,a)$$

$$\Big/(\Gamma(3+2\,n-2\,a+2\,b)\,\Gamma(\tfrac{3}{2}+b)\,\Gamma(-\tfrac{1}{2}+a-b-n)\,\Gamma(2\,a-1)\,\Gamma(a-b-1-n)\,\Gamma(-1-2\,b))$$

excess= $a$

B.16  Prudnikov 7.4.4.32  : T7

$${}_3F_2([a,\,b,\,n-b-1],\,[n-b-2+2\,a,\,b+1],\,1)\;=$$

$$-\frac{1}{2}\frac{\Gamma(b-n+2)\,\Gamma(a-b)\left(\sum_{L=1}^{n-2}\dfrac{\Gamma(-b+L)\,\Gamma(-a+b+2+L-n)}{\Gamma(2+b-n+L)\,\Gamma(-b+a+L)}\right)\Gamma(a)\,\Gamma(n-b-2+2\,a)}{\Gamma(2+b-a-n)\,\Gamma(n+a-b-1)\,\Gamma(-b)\,\Gamma(2\,a-1)}$$

$$-\frac{1}{2}\frac{\pi\,\Gamma(a-b)\,\Gamma(n-b-2+2\,a)\,\Gamma(b-n+2)}{\sin(\pi\,b)\,\Gamma(-b)\,\Gamma(n+2\,a-2\,b-2)\,\Gamma(a)\,\Gamma(-n+2\,b+2)}$$

excess= $a$

B.17  Prudnikov 7.4.4.22  : T7

$${}_3F_2([a,\,b,\,-b+2],\,[c,\,2\,a+2-c],\,1)\;=$$

$$\left(\left(\frac{1}{2}\frac{\Gamma(a-\tfrac{1}{2}c+1-\tfrac{1}{2}b)\,\Gamma(\tfrac{1}{2}b-1+\tfrac{1}{2}c)}{\Gamma(a+\tfrac{1}{2}b-\tfrac{1}{2}c)\,\Gamma(-\tfrac{1}{2}b+\tfrac{1}{2}c)}-\frac{1}{2}\frac{\Gamma(\tfrac{3}{2}+a-\tfrac{1}{2}b-\tfrac{1}{2}c)\,\Gamma(\tfrac{1}{2}b+\tfrac{1}{2}c-\tfrac{1}{2})}{\Gamma(a-\tfrac{1}{2}c+\tfrac{1}{2}+\tfrac{1}{2}b)\,\Gamma(-\tfrac{1}{2}b+\tfrac{1}{2}c+\tfrac{1}{2})}\right)\right.$$

$$\left.\Gamma(2\,a+2-c)\,\Gamma(c)\right)\Big/((-b+1)\,(-c+a+1)\,\Gamma(c+b-2)\,\Gamma(2+2\,a-b-c))$$

excess= $-\tfrac{1}{2}b+\tfrac{1}{2}n$

B.18  Prudnikov 7.4.4.32  : T1

$${}_3F_2([a,\,b,\,\tfrac{1}{2}a+\tfrac{1}{2}],\,[1+a,\,\tfrac{1}{2}n+\tfrac{1}{2}a-\tfrac{1}{2}+\tfrac{1}{2}b],\,1)\;=$$

$$\left( -\frac{1}{2}\cos(\frac{1}{2}\pi(b-a-n))\,\Gamma(\frac{3}{2}+\frac{1}{2}a-\frac{1}{2}n-\frac{1}{2}b)\,a\Gamma(\frac{1}{2}n+\frac{1}{2}a-\frac{1}{2}+\frac{1}{2}b)\sin(\frac{1}{2}\pi(b+n)) \right.$$

$$\times \left( \sum_{L=1}^{n-2} \frac{\Gamma(-\frac{1}{2}a-\frac{1}{2}n+\frac{1}{2}+\frac{1}{2}b+L)\,\Gamma(-\frac{1}{2}b-\frac{1}{2}n+1+L)}{\Gamma(\frac{3}{2}+\frac{1}{2}a-\frac{1}{2}n-\frac{1}{2}b+L)\,\Gamma(1+\frac{1}{2}b-\frac{1}{2}n+L)} \right) /\pi$$

$$\left. +\frac{2^{(a-1)}\sqrt{\pi}\,\Gamma(1+\frac{1}{2}a)\,\Gamma(\frac{1}{2}n+\frac{1}{2}a-\frac{1}{2}+\frac{1}{2}b)\,\Gamma(\frac{3}{2}+\frac{1}{2}a-\frac{1}{2}n-\frac{1}{2}b)}{\Gamma(1+a-b)\,\Gamma(\frac{1}{2}a+\frac{1}{2})\,\Gamma(b)} \right) \Big/ \sin(\frac{1}{2}\pi(-b+n))$$

$$\text{excess}= \frac{1}{2}n+b-\frac{1}{2}a-1$$

### B.19  Prudnikov 7.4.4.32  : T3

$${}_3F_2([a,\,b,\,1],\,[2\,b,\,\frac{1}{2}a+\frac{1}{2}n],\,1) \;=\;$$

$$\Gamma(1+b-\frac{1}{2}n-\frac{1}{2}a)\,\Gamma(\frac{1}{2}a-\frac{1}{2}n+1)$$

$$\times \left( \sum_{L=1}^{n-2} \frac{\Gamma(1-b-\frac{1}{2}n+\frac{1}{2}a+L)\,\Gamma(-\frac{1}{2}a-\frac{1}{2}n+1+L)}{\Gamma(1+b-\frac{1}{2}n-\frac{1}{2}a+L)\,\Gamma(1+\frac{1}{2}a-\frac{1}{2}n+L)} \right) (b-\frac{1}{2})$$

$$\Big/ (\Gamma(-b-\frac{1}{2}n+\frac{1}{2}a+2)\,\Gamma(-\frac{1}{2}n-\frac{1}{2}a+1))$$

$$+\frac{\sqrt{\pi}\,\Gamma(n-a+2\,b)\,\Gamma(a+n)\,\Gamma(b+\frac{1}{2})\,\Gamma(1+b-\frac{1}{2}n-\frac{1}{2}a)\,\Gamma(\frac{1}{2}a-\frac{1}{2}n+1)}{\Gamma(a)\,\Gamma(b)\,2^{(2\,n)}\,\Gamma(b+\frac{1}{2}n+\frac{1}{2}-\frac{1}{2}a)\,\Gamma(\frac{1}{2}a+\frac{1}{2}n+\frac{1}{2})\,(\frac{1}{2}n+b-\frac{1}{2}a-1)\,\Gamma(2\,b-a)}$$

$$\text{excess}= n-2\,b$$

### B.20  Prudnikov 7.4.4.32  : T6

$${}_3F_2([a,\,\frac{1}{2}n+\frac{1}{2}a-1,\,b],\,[n+a-b-1,\,\frac{1}{2}a+\frac{1}{2}n],\,1) \;=\;$$

$$\frac{2^{(-2-2b)}\,\Gamma(a+n)\,\Gamma(\tfrac{1}{2}n+\tfrac{1}{2}a-b-\tfrac{1}{2})\,\Gamma(-b+1)\,\Gamma(\tfrac{1}{2}a-\tfrac{1}{2}n+1)\,\Gamma(n+a-b-1)}{\Gamma(a)\,\Gamma(\tfrac{1}{2}n+\tfrac{1}{2}a-b)\,\Gamma(\tfrac{1}{2}a+\tfrac{1}{2}n+\tfrac{1}{2})\,\Gamma(n+a-2b-1)}$$

$$-\tfrac{1}{2}\Gamma(-b+1)\,\Gamma(\tfrac{1}{2}a-\tfrac{1}{2}n+1)\left(\sum_{L=1}^{n-2}\frac{\Gamma(b+1+L-n)\,\Gamma(-\tfrac{1}{2}a-\tfrac{1}{2}n+1+L)}{\Gamma(1-b+L)\,\Gamma(1+\tfrac{1}{2}a-\tfrac{1}{2}n+L)}\right)$$

$$\times \Gamma(n-2b)\,\Gamma(n+a-b-1)\Big/\big(\Gamma(-n+b+1)\,\Gamma(-\tfrac{1}{2}n-\tfrac{1}{2}a+1)\,\Gamma(n+a-2b-1)\,\Gamma(-b+n)\big)$$

$$\text{excess} = -\tfrac{1}{4}n - \tfrac{1}{2}a + 1$$

### B.21  Prudnikov 7.4.4.25 variation 2  : T2

$$_3F_2([2a, a, a+\tfrac{1}{2}], [\tfrac{1}{2}+\tfrac{3}{2}a-\tfrac{1}{4}n, 2a+1], 1) =$$

$$\frac{\left(\displaystyle\sum_{k=0}^{n-1}\frac{(-1)^k\,\Gamma(-a-\tfrac{1}{2}n+1+k)}{\Gamma(a-\tfrac{1}{2}n+1+k)}\right)\,\Gamma(\tfrac{1}{2}+\tfrac{3}{2}a-\tfrac{1}{4}n)\,2^{(2a)}}{\Gamma(-\tfrac{1}{2}a-\tfrac{1}{4}n+\tfrac{1}{2})}$$

$$+\frac{2^{(3a-1/2\,n-1)}\,\Gamma(a+\tfrac{1}{2}n)\,\Gamma(\tfrac{1}{2}+\tfrac{3}{2}a-\tfrac{1}{4}n)\,\Gamma(-\tfrac{1}{4}n-\tfrac{1}{2}a+1)}{\Gamma(2a)\,\sqrt{\pi}}$$

$$\text{excess} = 1 + \tfrac{1}{2}n + 2a$$

### B.22  Prudnikov 7.4.4.25 variation 2  : T3

$$_3F_2([a, 2a+\tfrac{1}{2}n, \tfrac{1}{2}-a-\tfrac{1}{2}n], [\tfrac{1}{2}+3a+\tfrac{1}{2}n, 1+a], 1) =$$

$$\tfrac{1}{2}\frac{\left(\displaystyle\sum_{k=0}^{n-1}\frac{(-1)^k\,\Gamma(1-2a-n+k)}{\Gamma(1+2a+k)}\right)\,\Gamma(\tfrac{1}{2}+3a+\tfrac{1}{2}n)\,\Gamma(2a+2)\,\Gamma(1+\tfrac{1}{2}n+2a)}{\Gamma(1-2a-n)\,\Gamma(1+4a+n)\,\Gamma(a+\tfrac{3}{2})}$$

$$+\frac{4^{(1/2\,n-1+3a)}\,\Gamma(2a+n)\,\Gamma(\tfrac{1}{2}+3a+\tfrac{1}{2}n)\,\Gamma(1+\tfrac{1}{2}n+2a)\,\Gamma(1+a)}{\sqrt{\pi}\,(2a+\tfrac{1}{2}n)\,\Gamma(4a+n)^2}$$

$$\text{excess} = -\frac{1}{2} + 2a + \frac{1}{2}n$$

**B.24 Prudnikov 7.4.4.25 variation 2 : T4**

$$_3F_2([a, -\frac{1}{2} + 2a + \frac{1}{2}n, -a - \frac{1}{2}n + 1], [\frac{1}{2}n - 1 + 3a, 1 + a], 1) =$$

$$\left( \frac{\left( \sum_{k=0}^{n-1} \frac{(-1)^k \Gamma(-2a - n + 2 + k)}{\Gamma(2a + k)} \right) 2^{(2a-2)}}{\Gamma(-2a - n + 2)\Gamma(-1 + 4a + n)} + \frac{4^{(1/2n + 3a - 3)} \Gamma(-1 + 2a + n)}{(-1 + 2a + \frac{1}{2}n)\Gamma(-2 + 4a + n)^2} \right) 2/\sqrt{\pi}$$

$$\times \Gamma(-\frac{1}{2} + 2a + \frac{1}{2}n)\Gamma(\frac{1}{2}n - 1 + 3a)\Gamma(1 + a)$$

$$\text{excess} = -\frac{1}{4}n - \frac{1}{4}a + \frac{1}{2}$$

**B.24 Prudnikov 7.4.4.25 variation 2 : T5**

$$_3F_2([a, \frac{1}{2}a + \frac{1}{2}, 1 + \frac{1}{2}a], [1 - \frac{1}{4}n + \frac{3}{4}a, 1 + a], 1) =$$

$$\frac{1}{2} \frac{\left( \sum_{k=0}^{n-1} \frac{(-1)^k \Gamma(-\frac{1}{2}n - \frac{1}{2}a + 1 + k)}{\Gamma(\frac{1}{2}a - \frac{1}{2}n + 1 + k)} \right) \Gamma(1 - \frac{1}{4}n + \frac{3}{4}a)\Gamma(\frac{1}{2}a - \frac{1}{2}n + 2)\Gamma(-\frac{1}{4}n - \frac{1}{4}a + \frac{1}{2})}{\Gamma(-\frac{1}{2}n - \frac{1}{2}a + 1)\Gamma(\frac{1}{4}a - \frac{1}{4}n + 1)\Gamma(\frac{3}{2} + \frac{1}{4}a - \frac{1}{4}n)}$$

$$+ \frac{2^{(3/2a - 1/2n - 1)} \Gamma(\frac{1}{2}a + \frac{1}{2}n)\Gamma(1 - \frac{1}{4}n + \frac{3}{4}a)\Gamma(-\frac{1}{4}n - \frac{1}{4}a + \frac{1}{2})}{\Gamma(a)\sqrt{\pi}}$$

$$\text{excess} = \frac{1}{2} + 2a + \frac{1}{2}n$$

**B.25 Prudnikov 7.4.4.25 variation 2 : T6**

$$_3F_2([a, \frac{1}{2} + 2a + \frac{1}{2}n, -a - \frac{1}{2}n + 1], [1 + 3a + \frac{1}{2}n, 1 + a], 1) =$$

$$\frac{1}{2} \frac{\left( \sum_{k=0}^{n-1} \frac{(-1)^k \Gamma(1 - 2a - n + k)}{\Gamma(1 + 2a + k)} \right) \Gamma(1 + 3a + \frac{1}{2}n)\Gamma(\frac{1}{2} + 2a + \frac{1}{2}n)\Gamma(2a + 2)}{\Gamma(1 - 2a - n)\Gamma(1 + 4a + n)\Gamma(a + \frac{3}{2})}$$

$$+ \frac{2^{(-2 + n + 6a)} \Gamma(2a + n)\Gamma(1 + 3a + \frac{1}{2}n)\Gamma(\frac{1}{2} + 2a + \frac{1}{2}n)\Gamma(1 + a)}{\sqrt{\pi}(2a + \frac{1}{2}n)\Gamma(4a + n)^2}$$

$$\text{excess} = -1 + 2a + \frac{1}{2}n$$

B.26  Prudnikov 7.4.4.25 variation 2   : T7

$$_3F_2([a, 2a + \frac{1}{2}n, -a + \frac{3}{2} - \frac{1}{2}n], [-\frac{1}{2} + 3a + \frac{1}{2}n, 1 + a], 1) =$$

$$\frac{1}{2} \frac{\left(\sum_{k=0}^{n-1} \frac{(-1)^k \Gamma(-2a - n + 2 + k)}{\Gamma(2a + k)}\right) \Gamma(-\frac{1}{2} + 3a + \frac{1}{2}n) \Gamma(-1 + 2a + \frac{1}{2}n) \Gamma(2a + 1)}{\Gamma(-2a - n + 2) \Gamma(-1 + 4a + n) \Gamma(a + \frac{1}{2})}$$

$$+ \frac{2^{(-5+n+6a)} \Gamma(1 + a) \Gamma(-1 + 2a + \frac{1}{2}n) \Gamma(-1 + 2a + n) \Gamma(-\frac{1}{2} + 3a + \frac{1}{2}n)}{(-1 + 2a + \frac{1}{2}n) \sqrt{\pi} \Gamma(-2 + 4a + n)^2}$$

$$\text{excess} = \frac{1}{2} + \frac{1}{2}a - \frac{1}{4}n$$

B.27  Prudnikov 7.4.4.25 variation 2   : T8

$$_3F_2([a, a + \frac{1}{2}, 1], [2a + 1, \frac{1}{2}a - \frac{1}{4}n + 1], 1) =$$

$$\frac{2^{(2a-2)} \Gamma(\frac{1}{2} + \frac{1}{2}a - \frac{1}{4}n) \Gamma(2 + a - \frac{1}{2}n) \Gamma(a + \frac{1}{2}n)}{\Gamma(2a) \Gamma(\frac{1}{2}a - \frac{1}{4}n + \frac{3}{2})}$$

$$+ \frac{\frac{1}{2}\left(\sum_{k=0}^{n-1} \frac{(-1)^k \Gamma(-a - \frac{1}{2}n + 1 + k)}{\Gamma(a - \frac{1}{2}n + 1 + k)}\right) \Gamma(\frac{1}{2} + \frac{1}{2}a - \frac{1}{4}n) \Gamma(2 + a - \frac{1}{2}n)}{\Gamma(-a - \frac{1}{2}n + 1) \Gamma(\frac{1}{2}a - \frac{1}{4}n + \frac{3}{2})}$$

$$\text{excess} = \frac{1}{2}a - \frac{1}{4}n - \frac{1}{2}$$

B.28  Prudnikov 7.4.4.25 variation 2   : T9

$$_3F_2([a - \frac{1}{2}, a, 1], [2a - 1, \frac{1}{2}a - \frac{1}{4}n + 1], 1) =$$

$$\left( \frac{2^{(-1+a-1/2\,n)} \left( \sum_{k=0}^{n-1} \frac{(-1)^k \, \Gamma(-\frac{1}{2}\,n-a+2+k)}{\Gamma(a-\frac{1}{2}\,n+k)} \right) \Gamma(\frac{1}{2}\,a-\frac{1}{4}\,n+1)}{\Gamma(-a-\frac{1}{2}\,n+2)\sqrt{\pi}} \right.$$

$$\left. + \frac{2^{(3\,a-4-1/2\,n)} \, \Gamma(-1+\frac{1}{2}\,n+a) \, \Gamma(\frac{1}{2}\,a-\frac{1}{4}\,n+1)}{\sqrt{\pi}\,\Gamma(2\,a-2)} \right) \Gamma(\frac{1}{2}\,a-\frac{1}{4}\,n-\frac{1}{2})$$

$$\text{excess} = 1 - \frac{1}{2}\,a + \frac{1}{2}\,n$$

### B.29  Prudnikov 7.4.4.25 variation 3  : T2

$${}_3F_2([2\,a-1,\,a,\,a-\frac{1}{2}],\,[-\frac{1}{2}+\frac{3}{2}\,a+\frac{1}{2}\,n,\,2\,a],\,1) =$$

$$-\frac{1}{2}\frac{\Gamma(1+a+n)\,\Gamma(a-n) \left( \sum_{k=1}^{2\,n-1} \frac{(-1)^k\,\Gamma(-a-n+1+k)}{\Gamma(a-n+k)} \right) \Gamma(-\frac{1}{2}+\frac{3}{2}\,a+\frac{1}{2}\,n)\,\Gamma(1-\frac{1}{2}\,a+\frac{1}{2}\,n)}{\Gamma(-a-n+1)\,\Gamma(a+n)\,\Gamma(\frac{1}{2}\,a+\frac{1}{2}\,n+1)\,\Gamma(\frac{1}{2}\,a+\frac{1}{2}\,n+\frac{1}{2})}$$

$$+ \frac{2^{(3\,a+n-3)}\,\Gamma(-\frac{1}{2}+\frac{3}{2}\,a+\frac{1}{2}\,n)\,\Gamma(a-n)\,\Gamma(1-\frac{1}{2}\,a+\frac{1}{2}\,n)}{\sqrt{\pi}\,\Gamma(2\,a-1)}$$

$$\text{excess} = -\frac{1}{2}\,n + 2\,a$$

### B.30  Prudnikov 7.4.4.25 variation 3  : T3

$${}_3F_2([a,\,-\frac{1}{2}\,n+2\,a,\,\frac{1}{2}+\frac{1}{2}\,n-a],\,[-\frac{1}{2}-\frac{1}{2}\,n+3\,a,\,1+a],\,1) =$$

$$-\frac{1}{2}\frac{\Gamma(2\,a+1)\,\Gamma(-n+2\,a) \left( \sum_{k=1}^{n-1} \frac{(-1)^k\,\Gamma(1-2\,a+k)}{\Gamma(2\,a-n+k)} \right) \Gamma(-\frac{1}{2}-\frac{1}{2}\,n+3\,a)\,\Gamma(-\frac{1}{2}\,n+2\,a)}{\Gamma(2\,a)\,\Gamma(a+\frac{1}{2})\,\Gamma(1-2\,a)\,\Gamma(4\,a-n)}$$

$$- \frac{2^{(6\,a-3-n)}\,(n-4\,a+1)\,\Gamma(1+a)\,\Gamma(-n+2\,a)\,\Gamma(-\frac{1}{2}-\frac{1}{2}\,n+3\,a)\,\Gamma(-\frac{1}{2}\,n+2\,a)}{\sqrt{\pi}\,\Gamma(4\,a-n)^2}$$

$$\text{excess} = \frac{3}{2} - \frac{1}{2}\,n + 2\,a$$

### B.31  Prudnikov 7.4.4.25 variation 3  : T4

$$_3F_2([a, -\tfrac{1}{2}n + \tfrac{1}{2} + 2a, -a + \tfrac{1}{2}n], [1 - \tfrac{1}{2}n + 3a, 1 + a], 1) =$$

$$-\frac{1}{2} \frac{\Gamma(2a+2)\,\Gamma(1+2a-n)\left(\sum_{k=1}^{n-1} \frac{(-1)^k \Gamma(-2a+k)}{\Gamma(1-n+2a+k)}\right) \Gamma(1 - \tfrac{1}{2}n + 3a)\,\Gamma(\tfrac{3}{2} - \tfrac{1}{2}n + 2a)}{\Gamma(-2a)\,\Gamma(2-n+4a)\,\Gamma(2a+1)\,\Gamma(a+\tfrac{3}{2})}$$

$$-\frac{2^{(-n+6a)}(n-1-4a)\,\Gamma(1+2a-n)\,\Gamma(1-\tfrac{1}{2}n+3a)\,\Gamma(1+a)\,\Gamma(\tfrac{3}{2} - \tfrac{1}{2}n + 2a)}{\sqrt{\pi}\,\Gamma(2-n+4a)^2}$$

$$\text{excess} = \tfrac{1}{2} + \tfrac{n}{4} - \tfrac{a}{2}$$

### B.32  Prudnikov 7.4.4.25 variation 3  : T5

$$_3F_2([2a-1, a+\tfrac{1}{2}, a], [\tfrac{3}{2}a + \tfrac{1}{4}n, 2a], 1) =$$

$$-\tfrac{1}{2}\Gamma(a - \tfrac{1}{2}n)\left(\sum_{k=1}^{n-1} \frac{(-1)^k \Gamma(-a - \tfrac{1}{2}n + 1 + k)}{\Gamma(a - \tfrac{1}{2}n + k)}\right)\Gamma(\tfrac{3}{2}a + \tfrac{1}{4}n)\,\Gamma(-\tfrac{1}{2}a + \tfrac{1}{4}n + \tfrac{1}{2})$$

$$\times \Gamma(1 + a + \tfrac{1}{2}n)\Big/\big(\Gamma(-a - \tfrac{1}{2}n + 1)\,\Gamma(1 + \tfrac{1}{2}a + \tfrac{1}{4}n)\,\Gamma(\tfrac{1}{2} + \tfrac{1}{2}a + \tfrac{1}{4}n)\,\Gamma(a + \tfrac{1}{2}n)\big)$$

$$+ \frac{2^{(n/2 + 3a - 3)}\,\Gamma(\tfrac{3}{2}a + \tfrac{1}{4}n)\,\Gamma(a - \tfrac{1}{2}n)\,\Gamma(-\tfrac{1}{2}a + \tfrac{1}{4}n + \tfrac{1}{2})}{\sqrt{\pi}\,\Gamma(2a-1)}$$

$$\text{excess} = -\tfrac{1}{2} - \tfrac{1}{2}n + 2a$$

### B.33  Prudnikov 7.4.4.25 variation 3  : T6

$$_3F_2([a, -\tfrac{1}{2}n + \tfrac{1}{2} + 2a, \tfrac{1}{2}n + 1 - a], [-\tfrac{1}{2}n + 3a, 1 + a], 1) =$$

$$-\frac{1}{2}\frac{\Gamma(2a+1)\,\Gamma(-n+2a)\left(\sum_{k=1}^{n-1} \frac{(-1)^k \Gamma(1-2a+k)}{\Gamma(2a-n+k)}\right)\Gamma(-\tfrac{1}{2} - \tfrac{1}{2}n + 2a)\,\Gamma(-\tfrac{1}{2}n + 3a)}{\Gamma(1-2a)\,\Gamma(4a-n)\,\Gamma(2a)\,\Gamma(a+\tfrac{1}{2})}$$

$$+ \frac{2^{(-2-n+6a)}\,\Gamma(-\tfrac{1}{2}n + \tfrac{1}{2} + 2a)\,\Gamma(1+a)\,\Gamma(-n+2a)\,\Gamma(-\tfrac{1}{2}n + 3a)}{\sqrt{\pi}\,\Gamma(4a-n)^2}$$

$$\text{excess} = -\frac{1}{2}n + 2a + 1$$

B.34  Prudnikov 7.4.4.25 variation 3  : T7

$$_3F_2([a, -\frac{1}{2}n + 2a + 1, \frac{1}{2} + \frac{1}{2}n - a], [-\frac{1}{2}n + 3a + \frac{3}{2}, 1 + a], 1) =$$

$$-\frac{1}{2} \frac{\Gamma(2a+2)\Gamma(1+2a-n)\left(\sum_{k=1}^{n-1} \frac{(-1)^k \Gamma(-2a+k)}{\Gamma(1-n+2a+k)}\right) \Gamma(-\frac{1}{2}n+2a+1)\Gamma(-\frac{1}{2}n+3a+\frac{3}{2})}{\Gamma(2a+1)\Gamma(a+\frac{3}{2})\Gamma(-2a)\Gamma(2-n+4a)}$$

$$-\frac{2^{(-n+6a)}(n-1-4a)\Gamma(1+a)\Gamma(1+2a-n)\Gamma(-\frac{1}{2}n+2a+1)\Gamma(-\frac{1}{2}n+3a+\frac{3}{2})}{\sqrt{\pi}\,\Gamma(2-n+4a)^2}$$

$$\text{excess} = \frac{1}{2}a + \frac{1}{4}n - \frac{1}{2}$$

B.35  Prudnikov 7.4.4.25 variation 3  : T8

$$_3F_2([a, a+\frac{1}{2}, 1], [2a, 1+\frac{1}{2}a+\frac{1}{4}n], 1) =$$

$$\left(-\frac{1}{2} \frac{\left(\sum_{k=1}^{n-1} \frac{(-1)^k \Gamma(-a-\frac{1}{2}n+1+k)}{\Gamma(a-\frac{1}{2}n+k)}\right) \Gamma(\frac{1}{2}a+\frac{1}{4}n-\frac{1}{2})}{\Gamma(-a-\frac{1}{2}n+1)\Gamma(a+\frac{1}{2}n)\Gamma(\frac{1}{2}+\frac{1}{2}a+\frac{1}{4}n)} + \frac{4^{(a-1)}}{(-1+\frac{1}{2}n+a)\Gamma(2a-1)}\right)$$

$$\Gamma(1+a+\frac{1}{2}n)\Gamma(a-\frac{1}{2}n)$$

$$\text{excess} = \frac{1}{2}a + \frac{1}{4}n$$

B.36  Prudnikov 7.4.4.25 variation 3  : T9

$$_3F_2([a-\frac{1}{2}, a, 1], [2a, \frac{1}{2}+\frac{1}{2}a+\frac{1}{4}n], 1) =$$

$$\left(-\frac{\left(\sum_{k=1}^{n-1} \frac{(-1)^k \Gamma(-a-\frac{1}{2}n+1+k)}{\Gamma(a-\frac{1}{2}n+k)}\right) \sin(\pi(a+\frac{1}{2}n))}{\pi} + \frac{4^{(a-1)}}{\Gamma(2a-1)}\right) \Gamma(a+\frac{1}{2}n)\Gamma(a-\frac{1}{2}n)$$

$$\text{excess} = -b + 2 + c - a$$

**B.37  Prudnikov 7.4.4.17  : T2**

$$_3F_2([a, 2, b], [c, 4], 1) =$$

$$-6\,\frac{(2c - 5 + b - ab + a)\,\Gamma(c)\,\Gamma(-b + 2 + c - a)}{(a - 3)(b - 1)(-2 + a)(a - 1)(b - 3)(b - 2)\,\Gamma(c - a)\,\Gamma(c - b)}$$

$$+\,\frac{6(c - 2)(c - 1)(ab - 3b - 3a + 3 + 2c)}{(a - 3)(b - 1)(-2 + a)(a - 1)(b - 3)(b - 2)}$$

$$\text{excess} = -b + 2 + c - a$$

**B.38  Gessel and Stanton, SIAM J. Math. Anal.,13,295(1982)  Eq.(5.16)  : T6**

$$_3F_2([a, a - n - 1, b], [c, a - n + 1], 1) =$$

$$-\Bigg((a - n)(a - n - 1)(-bn - ac + nc + a^2 - an + c - a + n)\,\Gamma(-b + 2 + c - a)\,\Gamma(c)\,\Gamma(n)$$

$$\left(\sum_{k=0}^{n}\frac{\Gamma(a - 1 - n + k)\,\Gamma(c - a + k)}{\Gamma(k + 1)\,\Gamma(1 - n + c - b + k)}\right)\Bigg)\Big/((b - 1)\,\Gamma(n + c + 1 - a)\,\Gamma(c - a)\,\Gamma(a))$$

$$-\,\frac{(a - n)(a - n - 1)\,\Gamma(-b + 2 + c - a)\,\Gamma(c)}{n(b - 1)\,\Gamma(c + 1 - b)\,\Gamma(c - a)}$$

$$\text{excess} = 2$$

**B.39  Prudnikov 7.4.4.17  : T3**

$$_3F_2([a, 2, b], [c, 4 + a - c + b], 1) =$$

$$\frac{(-a - 2 + c - b)(-a - 3 + c - b)(c - 2)(c - 1)(3 - 4c - ca + c^2 + 3a + ab + 3b - bc)}{(-a - 2 + c)(c - 3 - b)(c - 2 - b)(c - 1 - b)(-a - 3 + c)(-a + c - 1)}$$

$$+\,\frac{(a + 5 - ca + c^2 - 4c + ab - bc + b)\,\Gamma(4 + a - c + b)\,\Gamma(c)}{(-a - 2 + c)(c - 3 - b)(c - 2 - b)(c - 1 - b)(-a - 3 + c)(-a + c - 1)\,\Gamma(b)\,\Gamma(a)}$$

$$\text{excess} = 2$$

**B.40  Gessel & Stanton, SIAM J. Math. Anal. 13,2,295(1982)  Eq(1.9)  : T7**

$$_3F_2\left([a, b, -n + \frac{an}{b + n}], [1 + b + \frac{an}{b + n}, a - n + 1], 1\right) =$$

$$\frac{(-1)^n \sin(\pi b)\,\Gamma(1 + b + \frac{an}{b + n})\,\Gamma(-b)}{\sin(\pi a)\,\Gamma(1 + a)\,\Gamma(1 - n + \frac{an}{b + n})\,\Gamma(n - a)}$$

$$\text{excess} = 2$$

```
B.41  Gessel & Stanton, SIAM J. Math. Anal. 13,2,295(1982)  Eq(1.9)
with n->-1-n    : T7
```

$$_3F_2([a, b, -\frac{n(n-b+a)}{b-n}], [1+b-\frac{an}{b-n}, 1+a+n], 1) =$$

$$\frac{(-1)^n \sin(\pi b)\,\Gamma(1+b+\frac{an}{-b+n})\,\Gamma(-b)}{\sin(\pi a)\,\Gamma(1+a)\,\Gamma(-a-n)\,\Gamma(n+1+\frac{an}{-b+n})}$$

$$\text{excess} = 2$$

```
B.42  Gessel and Stanton, SIAM J. Math. Anal.,13,295(1982)
Eq.(5.16)    : T5}
```

$$_3F_2([a, b, c], [n+1+c+b, a-n+1], 1) =$$

$$\Gamma(n+1+c+b)\,\Gamma(a-n+1)\,(-n^2c+n^2a-n^3-n^2b-ncb)\,\Gamma(n)$$
$$\left(\sum_{k=0}^{n}\frac{\Gamma(b+k)\,\Gamma(c+k)}{\Gamma(k+1)\,\Gamma(a+1+k-n)}\right)$$
$$/((a-c-b-n)\,\Gamma(c)\,\Gamma(b)\,\Gamma(c+n+1)\,\Gamma(b+n+1)\,n)$$
$$+\frac{\Gamma(n+1+c+b)\,\Gamma(a-n+1)}{n\,(a-c-b-n)\,\Gamma(1+a)\,\Gamma(c)\,\Gamma(b)}$$

$$\text{excess} = c - \frac{1}{2}a - \frac{1}{2}b - \frac{1}{2}$$

```
B.43  Prudnikov 7.4.4.20  : T1
```

$$_3F_2([a, b, c], [2c-1, \frac{1}{2}a+\frac{1}{2}b+\frac{1}{2}], 1) =$$

$$\sqrt{\pi}\,\Gamma(c-\frac{1}{2})\,\Gamma(\frac{1}{2}a+\frac{1}{2}b+\frac{1}{2})\,\Gamma(c-\frac{1}{2}a-\frac{1}{2}b-\frac{1}{2})(\Gamma(\frac{1}{2}a)\,\Gamma(\frac{1}{2}b)\,\Gamma(c-\frac{1}{2}a)\,\Gamma(c-\frac{1}{2}b)$$
$$+\Gamma(\frac{1}{2}b+\frac{1}{2})\,\Gamma(c-\frac{1}{2}-\frac{1}{2}a)\,\Gamma(c-\frac{1}{2}-\frac{1}{2}b)\,\Gamma(\frac{1}{2}a+\frac{1}{2}))$$
$$\bigg/(\Gamma(\frac{1}{2}b+\frac{1}{2})\Gamma(c-\frac{1}{2}-\frac{1}{2}a)\,\Gamma(c-\frac{1}{2}-\frac{1}{2}b)\,\Gamma(\frac{1}{2}a+\frac{1}{2})\,\Gamma(c-\frac{1}{2}a)\,\Gamma(c-\frac{1}{2}b)\,\Gamma(\frac{1}{2}a)\,\Gamma(\frac{1}{2}b))$$

$$\text{excess} = -\frac{1}{2}c - \frac{1}{2}a + 1 + b$$

```
B.44  Prudnikov 7.4.4.20  : T2
```

$$_3F_2([a, b, c], [1+2b, \frac{1}{2}c+\frac{1}{2}a], 1) =$$

$$\left(\left(\frac{\Gamma(\frac{1}{2}a)}{\Gamma(1+b-\frac{1}{2}c)\,\Gamma(-\frac{1}{2}a+b+\frac{1}{2})\,\Gamma(\frac{1}{2}c)} + \frac{\Gamma(\frac{1}{2}a+\frac{1}{2})}{\Gamma(\frac{1}{2}+b-\frac{1}{2}c)\,\Gamma(-\frac{1}{2}a+1+b)\,\Gamma(\frac{1}{2}c+\frac{1}{2})}\right)\right.$$

$$\left.\times \sqrt{\pi}\,\Gamma(\frac{1}{2}c+\frac{1}{2}a)\,\Gamma(1+2b)\,\Gamma(-\frac{1}{2}c-\frac{1}{2}a+1+b)\right)\Big/(2^{(-a+1+2b)}\,\Gamma(b+1)\,\Gamma(a))$$

$$\text{excess}= a - 1$$

### B.45  Prudnikov 7.4.4.20  : T7

$$_3F_2([a, b, -b+1], [c, 2a-c], 1) =$$

$$\left(\frac{\sqrt{\pi}\,\Gamma(a-1)\,\Gamma(c)\,\Gamma(\frac{1}{2}b+\frac{1}{2}c-\frac{1}{2})}{2^{(2a-b-c)}\,\Gamma(a)\,\Gamma(-\frac{1}{2}b+\frac{1}{2}c+\frac{1}{2})\,\Gamma(a-\frac{1}{2}b-\frac{1}{2}c)\,\Gamma(a+\frac{1}{2}b-\frac{1}{2}c-\frac{1}{2})\,\Gamma(c+b-1)}\right.$$

$$\left.+ \frac{\sqrt{\pi}\,\Gamma(\frac{1}{2}b+\frac{1}{2}c)\,\Gamma(a-1)\,\Gamma(c)}{2^{(2a-b-c)}\,\Gamma(a)\,\Gamma(-\frac{1}{2}b+\frac{1}{2}c)\,\Gamma(-\frac{1}{2}b+a-\frac{1}{2}c+\frac{1}{2})\,\Gamma(a+\frac{1}{2}b-\frac{1}{2}c)\,\Gamma(c+b-1)}\right)\Gamma(2a-c)$$

$$\text{excess}= 1 + a$$

### B.46  Prudnikov 7.4.4.20  : T8

$$_3F_2([a, b, -b], [c, 2a-c+1], 1) =$$

$$\frac{\sqrt{\pi}\,\Gamma(2a-c+1)\,\Gamma(c)\,\Gamma(\frac{1}{2}b+\frac{1}{2}c)}{2^{(2a-c+1-b)}\,\Gamma(a-\frac{1}{2}c+1+\frac{1}{2}b)\,\Gamma(-\frac{1}{2}b+a-\frac{1}{2}c+\frac{1}{2})\,\Gamma(-\frac{1}{2}b+\frac{1}{2}c)\,\Gamma(c+b)}$$

$$+\frac{\sqrt{\pi}\,\Gamma(2a-c+1)\,\Gamma(\frac{1}{2}c+\frac{1}{2}b+\frac{1}{2})\,\Gamma(c)}{2^{(2a-c+1-b)}\,\Gamma(a-\frac{1}{2}c+\frac{1}{2}+\frac{1}{2}b)\,\Gamma(a-\frac{1}{2}c+1-\frac{1}{2}b)\,\Gamma(-\frac{1}{2}b+\frac{1}{2}c+\frac{1}{2})\,\Gamma(c+b)}$$

$$\text{excess}= -\frac{1}{2}a + 1 - \frac{1}{2}b + c$$

### B.47  Prudnikov 7.4.4.22  : T1

$$_3F_2([a, b, c], [2c, 1+\frac{1}{2}b+\frac{1}{2}a], 1) =$$

$$\left(\left(\frac{1}{2}\frac{\Gamma(c-\frac{1}{2}a)\Gamma(\frac{1}{2}a)}{\Gamma(c-\frac{1}{2}b)\Gamma(\frac{1}{2}b)} - \frac{1}{2}\frac{\Gamma(\frac{1}{2}-\frac{1}{2}a+c)\Gamma(\frac{1}{2}a+\frac{1}{2})}{\Gamma(\frac{1}{2}-\frac{1}{2}b+c)\Gamma(\frac{1}{2}b+\frac{1}{2})}\right)\Gamma(1+\frac{1}{2}b+\frac{1}{2}a)\Gamma(2c)\Gamma(-\frac{1}{2}a+1-\frac{1}{2}b+c)\right)$$
$$\Big/ \left(\Gamma(2c-a)\Gamma(c)\Gamma(a)(\frac{1}{2}b-\frac{1}{2}a)(-\frac{1}{2}a-\frac{1}{2}b+c)\right)$$

$$\text{excess} = -\frac{1}{2}c - \frac{1}{2}b + a - \frac{1}{2}n + \frac{1}{2} + \frac{1}{2}m$$

B.48  Lewanowicz, J. Comp & Appl. Math. 86,375(1997) Generalized
Watson; Eq.(2.15)  : T3

$$_3F_2([a, b, c], [\frac{1}{2}c + \frac{1}{2}b + \frac{1}{2}n + \frac{1}{2} + \frac{1}{2}m, -n+2a], 1) =$$

$$\frac{W_{m,n}(c, -b+2a-n, a-n, m, n)\Gamma(-\frac{1}{2}c-\frac{1}{2}b+a-\frac{1}{2}n+\frac{1}{2}+\frac{1}{2}m)\Gamma(\frac{1}{2}c+\frac{1}{2}b+\frac{1}{2}n+\frac{1}{2}+\frac{1}{2}m)}{\Gamma(-\frac{1}{2}c+\frac{1}{2}b+\frac{1}{2}n+\frac{1}{2}+\frac{1}{2}m)\Gamma(\frac{1}{2}c-\frac{1}{2}b+a-\frac{1}{2}n+\frac{1}{2}+\frac{1}{2}m)}$$

$$\text{excess} = b + n$$

B.49  Lewanowicz, J. Comp & Appl. Math. 86,375(1997) Generalized
Watson; Eq.(2.15)  : T9

$$_3F_2([a, -a+m+1, b], [c, -c+n+2b+1+m], 1) =$$

$$\frac{W_{m,n}(a-c+n+2b, -c-a+n+1+m+2b, b, m, n)\Gamma(b+n)\Gamma(c)}{\Gamma(c-b)\Gamma(2b+n)}$$

$$\text{excess} = b - 2c - 2a + 3$$

B.50  Lavoie, Math. Comp., 49,179,269(1987),  Eq(2)  : T2

$$_3F_2([a, b, c], [b+2-a, 1+b-c], 1) =$$

$$-\frac{\Gamma(b+2-a)\Gamma(1+b-c)\Gamma(\frac{1}{2}b-c-a+2)\sqrt{\pi}}{\Gamma(-a+1+\frac{1}{2}b)\Gamma(\frac{1}{2}b-c+1)\Gamma(\frac{1}{2}b+\frac{1}{2})(a-1)\Gamma(-c-a+b+2)2^b}$$

$$+\frac{\Gamma(b+2-a)\Gamma(1+b-c)\Gamma(\frac{1}{2}b-c-a+\frac{3}{2})\sqrt{\pi}}{(a-1)\Gamma(-c-a+b+2)\Gamma(\frac{1}{2}b)\Gamma(\frac{1}{2}b+\frac{3}{2}-a)2^b\Gamma(\frac{1}{2}+\frac{1}{2}b-c)}$$

$$\text{excess} = a - \frac{1}{2}b - \frac{1}{2}c$$

**B.51  Lavoie, Math. Comp., 49,179,269(1987), Eq(2)  : T3**

$$_3F_2([a, b, c], [1 + \tfrac{1}{2}c + \tfrac{1}{2}b, 2a - 1], 1) =$$

$$\left(\left(-\frac{1}{4}\frac{c(-b-1+2a)\Gamma(a-\tfrac{1}{2}b)\Gamma(\tfrac{1}{2}c+\tfrac{1}{2})}{\Gamma(-\tfrac{1}{2}c+a)\Gamma(\tfrac{1}{2}b+\tfrac{1}{2})} + \frac{\Gamma(1+\tfrac{1}{2}c)\Gamma(-\tfrac{1}{2}b+a+\tfrac{1}{2})}{\Gamma(-\tfrac{1}{2}c+a-\tfrac{1}{2})\Gamma(\tfrac{1}{2}b)}\right)\right.$$

$$\left.\times \Gamma(1+\tfrac{1}{2}c+\tfrac{1}{2}b)\Gamma(a-\tfrac{1}{2}b-\tfrac{1}{2}c)\Gamma(a-\tfrac{1}{2})\, 2^{(c-b+2a-1)}\right) \Big/ (\sqrt{\pi}\,\Gamma(c+1)\Gamma(2a-b)(-\tfrac{1}{2}c+\tfrac{1}{2}b))$$

$$\text{excess} = b + 1$$

**B.52  Lavoie, Math. Comp., 49,179,269(1987), Eq(2)  : T9**

$$_3F_2([a, -a+1, b], [c, -c+2+2b], 1) =$$

$$\left(-\frac{1}{(c-1-b)\Gamma(\tfrac{1}{2}+\tfrac{1}{2}a-\tfrac{1}{2}c+b)\Gamma(-\tfrac{1}{2}c-\tfrac{1}{2}a+1+b)\Gamma(\tfrac{1}{2}c+\tfrac{1}{2}a)\Gamma(-\tfrac{1}{2}a+\tfrac{1}{2}c+\tfrac{1}{2})}\right.$$

$$\left.+\frac{1}{(c-1-b)\Gamma(\tfrac{1}{2}a-\tfrac{1}{2}c+1+b)\Gamma(-\tfrac{1}{2}a+\tfrac{1}{2}c)\Gamma(-\tfrac{1}{2}+\tfrac{1}{2}c+\tfrac{1}{2}a)\Gamma(\tfrac{3}{2}-\tfrac{1}{2}c-\tfrac{1}{2}a+b)}\right)$$

$$\times \pi\, 2^{(-2b)}\,\Gamma(-c+2+2b)\,\Gamma(c)$$

$$\text{excess} = \frac{n+2-2a+ab-b}{n+1-a}$$

**B.53  Gessel & Stanton, SIAM J. Math. Anal. 13,2,295(1982)   Eq(1.9)**
    **: T5**

$$_3F_2([a, b, a-n-1], [a - \frac{n(b-1)}{a-n-1}, a-n+1], 1) =$$

$$\sin(\pi(b-a+\frac{n(-b+1)}{n+1-a}))(n-a)(bn-a+1)\Gamma(a+\frac{n(b-1)}{n+1-a})$$

$$\times \Gamma(2-a+\frac{(-a+1)(b-1)}{n+1-a})\Gamma(2-b+\frac{n(b-1)}{n+1-a}) \Big/ ((b-a+n)(a-1)\pi\,\Gamma(2+\frac{n(b-1)}{n+1-a}))$$

$$\text{excess} = -\frac{(b+n-2)(1+a-b)}{b-2}$$

**B.54  Gessel & Stanton, SIAM J. Math. Anal. 13,2,295(1982)   Eq(1.9)**
    **: T8**

$$_3F_2([2, a, 1-n], [b, 2 + \frac{n(-a+1)}{b-2}], 1) =$$

$$\frac{(a\,n - b - n + 2)\,(b - 1)}{(1 + a - b)\,(b + n - 2)}$$

excess= $n + 1$

```
B.55  Gessel & Stanton, SIAM J. Math. Anal. 13,2,295(1982)  Eq(1.9)
with n->-1-n   : T1
```

$${}_3\mathrm{F}_2([a, b, -\frac{(a-1)(b-1)}{n}], [2 + \frac{(b-1)(-a+1)}{n}, -1 + a + b + n], 1) =$$

$$\frac{\sin(\pi a)\,(-1)^n\,(n - ab + a + b - 1)\,\Gamma(n)\,\Gamma(-a - n + 1)\,\Gamma(-1 + a + b + n)}{\pi\,\Gamma(b + n)}$$

excess= $\dfrac{n\,(n - b + a)}{a + n - 1}$

```
B.56  Gessel & Stanton, SIAM J. Math. Anal. 13,2,295(1982)  Eq(1.9)
with n->-1-n   : T2
```

$${}_3\mathrm{F}_2([a, 2, b], [1 + \frac{b\,(a-1) + n}{a + n - 1}, 1 + a + n], 1) =$$

$$\frac{(n + a\,b - b)\,(a + n)}{(n - b + a)\,n}$$

excess= $-\dfrac{3\,b + a\,b - b^2 - 2\,a - 2 - n + b\,n}{1 + a - b}$

```
B.57  Gessel & Stanton, SIAM J. Math. Anal. 13,2,295(1982)  Eq(1.9)
with n->-1-n   : T4
```

$${}_3\mathrm{F}_2([a, b, 1 + a + \frac{a\,n}{1 + a - b}], [2 + a, 1 + a + n], 1) =$$

$$-\frac{\sin(\frac{\pi n a}{1 + a - b})\,(1 + a)\,\Gamma(2 - b - \frac{n\,(-b+1)}{b - a - 1})\,\Gamma(\frac{n\,(b-1)}{1 + a - b})}{\sin(\pi a)\,\Gamma(-b + a + 2 + n)\,\Gamma(-a - n)}$$

excess= $\dfrac{2\,a + n - 2 - a\,b + b}{a + n - 1}$

```
B.58  Gessel & Stanton, SIAM J. Math. Anal. 13,2,295(1982)  Eq(1.9)
with n->-1-n   : T5
```

$${}_3\mathrm{F}_2([a, b, a + n - 1], [a - 1 + \frac{b\,n + a - 1}{a + n - 1}, 1 + a + n], 1) =$$

$$-\sin(\pi\,(a + \frac{(-b+1)(a-1)}{a + n - 1}))\,(b\,n + a - 1)\,(a + n)\,\Gamma(1 + \frac{(b-1)(-a+1)}{a + n - 1})$$
$$\times \Gamma(a + \frac{n\,(b-1)}{a + n - 1})\,\Gamma(2 - a + \frac{(b-1)(a-1)}{a + n - 1}) \Big/ ((n - b + a)\,\pi\,(a - 1)\,\Gamma(\frac{b\,n + a - 1}{a + n - 1} + 1))$$

$$\text{excess} = \frac{(n+2-b)(1+a-b)}{b-2}$$

```
B.59  Gessel & Stanton, SIAM J. Math. Anal. 13,2,295(1982)  Eq(1.9)
with n->-1-n   : T8
```

$$_3F_2([2, a, n+1], [b, 2 + \frac{n(a-1)}{b-2}], 1) =$$

$$-\frac{(-1)^n \sin(\pi(b-a+\frac{n(a-1)}{b-2}))(an+b-2-n)(b-1)}{(1+a-b)(n+2-b)\sin(\frac{\pi(n+2-b)(1+a-b)}{b-2})}$$

```
The following are special cases of Minton type (Refs. 9 and 10)
that reduce an infinite series to a finite one.  See also B.42 and B.53.
```

$$\text{excess} = a - b - c - n + 3$$

```
B.60  Gessel and Stanton, SIAM J. Math. Anal.,13,295(1982)
Eq.(5.16)}  : T2
```

$$_3F_2([a, b, c], [2+a, a-n+1], 1) =$$

$$\left(\left(ncb - nc - anc - anb + n + 2an - bn - n^2a + na^2\right)\Gamma(2+a)\Gamma(a+2-b-c-n)\right.$$

$$\left.\Gamma(a-n+1)\Gamma(n)\left(\sum_{k=0}^{n}\frac{\Gamma(a+k+1-n-c)\Gamma(-b-n+1+a+k)}{\Gamma(k+1)\Gamma(a+1+k-n)}\right)\right)$$

$$/(n\,\Gamma(-b-n+1+a)\,\Gamma(2+a-b)\,\Gamma(2+a-c)\,\Gamma(1+a-c-n))$$

$$+ \frac{\Gamma(a+2-b-c-n)(-1-a)\Gamma(a-n+1)}{\Gamma(-b-n+1+a)\,n\,\Gamma(1+a-c-n)}$$

$$\text{excess}= c - 1 - n - b$$

**B.61  Gessel and Stanton, SIAM J. Math. Anal.,13,295(1982)**
`Eq.(5.16)    : T3`

$$_3F_2([a, b, 2], [c, a - n + 1], 1) =$$

$$-\Bigg((a - n - 1)(a - n)(a b - b - 1 - 2n + a - a c + n c + c)\Gamma(c)\Gamma(n)$$

$$\left(\sum_{k=0}^{n}\frac{\Gamma(a - 1 - n + k)\Gamma(c - 1 - n - b + k)}{\Gamma(k+1)\Gamma(c - n - 1 + k)}\right)\Bigg)/((a - n - b)\Gamma(a)\Gamma(c - b))$$

$$-\frac{(a - n - 1)(a - n)(c - 1)}{n(a - n - b)}$$

$$\text{excess}= 1$$

**B.62  Prudnikov 7.4.4.23  : T6**

$$_3F_2([a, b, c], [a - n + 1, c + b + n], 1) =$$

$$-\frac{\sin(\pi c)\left(\sum_{L=0}^{n-1}\frac{\Gamma(1 - c - n + L)\Gamma(b + L)}{\Gamma(L+1)\Gamma(b - a + 1 + L)}\right)\Gamma(c + b + n)\Gamma(b - a + n)\Gamma(a - n + 1)\Gamma(n)}{\pi\,\Gamma(b)\Gamma(a)\Gamma(b + n)}$$

$$\text{excess}= b - c - n + 2 - a$$

**B.63  Prudnikov 7.4.4.23  : T1**

$$_3F_2([a, b, c], [b + 1, -n + b + 1], 1) =$$

$$-(-1)^n b\left(\sum_{L=0}^{n-1}\frac{\Gamma(-b + a + L)\Gamma(-c - n + b + 1 + L)}{\Gamma(L+1)\Gamma(-c - n + 2 + L)}\right)\Gamma(b - c - n + 2 - a)\Gamma(-n + b + 1)$$

$$\Gamma(n)\Gamma(-c + 1)/(\Gamma(a - b)\Gamma(b - a + 1)\Gamma(b - c - n + 1)\Gamma(1 + b - c))$$

$$\text{excess}= c - a - n$$

**B.64  Prudnikov 7.4.4.23  : T2**

$$_3F_2([a, b, 1], [c, -n + b + 1], 1) =$$

$$\frac{(c - 1)\left(\sum_{L=0}^{n-1}\frac{\Gamma(a + L + 1 - c)\Gamma(-n + b + L)}{\Gamma(L+1)\Gamma(b - c - n + 2 + L)}\right)(-b + n)\Gamma(n)\Gamma(1 + b - c)}{\Gamma(n + 1 - c + a)\Gamma(b)}$$

$$\text{excess}= -\frac{n(b - a + n)}{n + 1 - a}$$

**B.65  Gessel & Stanton, SIAM J. Math. Anal. 13,2,295(1982)  Eq(1.9)**
`  : T2`

$$_3F_2([a, 2, b], [2 + \frac{(b-1)(a-1)}{a-n-1}, a-n+1], 1) =$$
$$\frac{(n-a)(n-ab+b)}{(b-a+n)n}$$

$$\text{excess} = \frac{-3b - ab + b^2 + 2a + 2 - n + bn}{1 + a - b}$$

```
B.66  Gessel & Stanton, SIAM J. Math. Anal. 13,2,295(1982)  Eq(1.9)
   : T4
```

$$_3F_2([a, b, 1 + a - \frac{an}{1+a-b}], [2 + a, a - n + 1], 1) =$$
$$\frac{(1+a)\sin(\frac{\pi a n}{1+a-b})\Gamma(2-b+\frac{n(b-1)}{1+a-b})\Gamma(-\frac{n(b-1)}{1+a-b})}{\sin(\pi a)\Gamma(-b+a+2-n)\Gamma(n-a)}$$

### Appendix C

The results contained in this Appendix together with the contents of Section 3.4 of the text are sufficient to allow the computation of any hypergeometric sum contiguous to any sum defined by Watson's, Whipple's or Dixon's theorems.

The following is the recursion in m for fixed n=N for the elements of the generalized Watson's theorem

$$W(m, n = N) =$$
$$-\frac{(m-1+a+b)(-3+a+b+m)(a+b-2c+3-m-2N)W(m-4, N)}{(-1+m+a+b-2c)(m-1-b+a)(-m+a-b+1)}$$
$$-((m-1+a+b)(10 + 2m^2 + 2mN + 4ca + 4cb - 2b^2 - 2a^2 + 2Na + 2Nb - 8m - 6N - 4c)$$
$$\times W(m-2, N))/((-1+m+a+b-2c)(m-1-b+a)(-m+a-b+1))$$
$$\tag{C.1}$$

The following is the recursion in n for fixed m=M for the generalized Watson's theorem.

$$W(m = M, n) = \frac{1}{2}((2c+n-1)(-3na - 3nb - 11n + 2Mc - 4ca - 4cb - 16c + 8$$
$$+ 12cn + 4a + 4b - 2M + 8c^2 + 4n^2 + 2ab + Mn)W(M, n-1))$$
$$/((-1+n+c)(-n+1-2c+b)(-n+1-2c+a))$$
$$+ \frac{\frac{1}{2}(2c+n-1)(2c+n-2)(a+b-2c+3-M-2n)W(M, n-2)}{(-1+n+c)(-n+1-2c+b)(-n+1-2c+a)}$$

(C.2)

This is the case m=-1, n=0 for Dixon's theorem, obtained by recursion as discussed in the text.

$$X(-1, 0) = \frac{2^{(-a)} \sqrt{\pi}\, \Gamma(a-c)\, \Gamma(a-b)\, \Gamma(\frac{1}{2}a - c - b)}{\Gamma(\frac{1}{2}a + \frac{1}{2})\, \Gamma(\frac{1}{2}a - b)\, \Gamma(\frac{1}{2}a - c)\, \Gamma(a - b - c)}$$
$$+ \frac{2^{(-a)} \sqrt{\pi}\, \Gamma(a-c)\, \Gamma(a-b)\, \Gamma(\frac{1}{2}a - c - b + \frac{1}{2})}{\Gamma(\frac{1}{2}a)\, \Gamma(\frac{1}{2} + \frac{1}{2}a - b)\, \Gamma(a - b - c)\, \Gamma(\frac{1}{2} + \frac{1}{2}a - c)}$$

(C.3)

This is the case m=2, n=0 for Dixon's theorem, obtained by recursion as discussed in the text.

$$X(2, 0) = -2\, \frac{\sqrt{\pi}\, \Gamma(\frac{7}{2} + \frac{1}{2}a - b - c)\, \Gamma(3 + a - c)\, \Gamma(3 + a - b)}{(-2 + b)(-2 + c)\, 2^a\, (b-1)(c-1)\, \Gamma(a + 3 - b - c)\, \Gamma(\frac{3}{2} + \frac{1}{2}a - c)\, \Gamma(\frac{3}{2} + \frac{1}{2}a - b)\, \Gamma(\frac{1}{2}a)}$$
$$+ ((\frac{1}{4}\, \frac{(2 + a - 2c)(-2b + a + 2)(a + 5 - 2c - 2b)\, \sqrt{\pi}}{(-2 + b)(-2 + c)\, 2^a\, (b-1)(c-1)}$$
$$- \frac{1}{8}\, \frac{(-2b - 2c + a + 6)\, \sqrt{\pi}}{(-2 + b)(-2 + c)\, 2^{(-2+a)}})\Gamma(2 + \frac{1}{2}a - c - b)\, \Gamma(3 + a - b)\, \Gamma(3 + a - c))$$
$$\Big/ (\Gamma(\frac{1}{2}a + \frac{1}{2})\, \Gamma(a + 3 - b - c)\, \Gamma(\frac{1}{2}a + 2 - c)\, \Gamma(\frac{1}{2}a - b + 2))$$

(C.4)

Eq.(21) then gives

$$W(-1, 0) = \left( \frac{1}{\Gamma(\frac{1}{2} - \frac{1}{2}a + c)\, \Gamma(c - \frac{1}{2}b)\, \Gamma(\frac{1}{2}b)\, \Gamma(\frac{1}{2}a + \frac{1}{2})} \right.$$
$$\left. + \frac{1}{\Gamma(\frac{1}{2} - \frac{1}{2}b + c)\, \Gamma(\frac{1}{2}b + \frac{1}{2})\, \Gamma(\frac{1}{2}a)\, \Gamma(c - \frac{1}{2}a)} \right) \Gamma(-\frac{1}{2}a - \frac{1}{2}b + c)\, \Gamma(\frac{1}{2}a + \frac{1}{2}b)\, \Gamma(c + \frac{1}{2})\, \sqrt{\pi}$$

(C.5)

$$\mathrm{W}(2,\,0) = \left(\left(-\frac{a^2 - 2\,c\,a - 2\,c\,b - 1 + 2\,c + b^2}{\Gamma(\frac{1}{2}\,a + \frac{1}{2})\,\Gamma(\frac{1}{2} - \frac{1}{2}\,a + c)\,\Gamma(\frac{1}{2} - \frac{1}{2}\,b + c)\,\Gamma(\frac{1}{2}\,b + \frac{1}{2})}\right.\right.$$

$$\left.\left. - \frac{8}{\Gamma(\frac{1}{2}\,b)\,\Gamma(\frac{1}{2}\,a)\,\Gamma(c - \frac{1}{2}\,b)\,\Gamma(c - \frac{1}{2}\,a)}\right)\sqrt{\pi}\,\Gamma(-\frac{1}{2} - \frac{1}{2}\,a - \frac{1}{2}\,b + c)\,\Gamma(\frac{1}{2}\,a + \frac{1}{2}\,b + \frac{3}{2})\right.$$

$$\left. \times \Gamma(c + \frac{1}{2}) \right) / ((-1 - a + b)\,(-a + b + 1))$$

(C.6)